\documentclass[12pt]{article}
\usepackage{xy}
\usepackage{amsfonts,amsmath}
\usepackage{amssymb,latexsym}
\usepackage[dvips]{epsfig}
\usepackage{amscd}
\usepackage{psfrag}

\usepackage[hmargin=3cm,vmargin=3.5cm]{geometry}

\title{Heisenberg algebra and a graphical calculus}
\author{Mikhail Khovanov}
\date{September 16, 2010}
 
\newtheorem{prop}{Proposition}
\newtheorem{theorem}{Theorem}
\newtheorem{cor}{Corollary}
\newtheorem{conj}{Conjecture}
\newtheorem{lemma}{Lemma}
\newcommand{\oplusop}[1]{{\mathop{\oplus}\limits_{#1}}}
\newcommand{\oplusoop}[2]{{\mathop{\oplus}\limits_{#1}^{#2}}}
 
\begin{document} 
     
\maketitle
\baselineskip 14pt
 
\tableofcontents

\def\R{\mathbb R}
\def\Q{\mathbb Q}
\def\Z{\mathbb Z}
\def\N{\mathbb N} 
\def\C{\mathbb C}
\def\l{\lbrace}
\def\r{\rbrace}
\def\o{\otimes}
\def\lra{\longrightarrow}
\def\RHom{\mathrm{RHom}}
\def\Id{\mathrm{Id}}
\def\mc{\mathcal}
\def\mf{\mathfrak} 
\def\Ext{\mathrm{Ext}}
\def\End{\mathrm{End}}
\def\Sym{\mathrm{Sym}}
\newcommand{\hei}{H} 
\newcommand{\Hei}{\mc{C}}
 \def\yesnocases#1#2#3#4{\left\{
\begin{array}{ll} #1 & #2 \\ #3 & #4
\end{array} \right. }

\def\gdim{{\mathrm{gdim}}}
\def\dmod{{\mathrm{-mod}}} 
\newcommand{\Hom}{{\rm Hom}}

\newcommand{\Ind}{{\rm Ind}}
\newcommand{\Res}{{\rm Res}}
\newcommand{\finone}{BFin}
\newcommand{\fintwo}{IRFin'}

\newcommand{\Cone}{\mc{H}'} %% basic category
\newcommand{\Ctwo}{\mc{H}} %% its Karoubian envelope

\newcommand{\Sone}{\mc{S}'} %%
\newcommand{\Stwo}{\mc{S}} %% its Karoubian envelope

\vspace{0.7in} 

\begin{abstract} 
A new calculus of planar diagrams involving diagrammatics 
for biadjoint functors and degenerate affine Hecke algebras is introduced. 
The calculus leads to an additive monoidal category whose Grothendieck 
ring contains an integral form of the Heisenberg algebra in infinitely 
many variables. We construct bases of vector spaces of morphisms 
between products of generating objects in this category.  
\end{abstract}

%%%%%%%%%%%%% PSFRAG%%%%%%%%%%%%%%%%%%%%
\psfrag{Qpm}{$Q_{+-}$} \psfrag{Qmp}{$Q_{-+}$}
\psfrag{Ione}{$\quad\mathbf{1}$}
\psfrag{Xiota}{$\iota$}\psfrag{Xiotap}{$\iota'$}
\psfrag{Xgi}{$g_i$}\psfrag{Xgim1}{$g_i^{-1}$}
\psfrag{Xsumi}{\large{$\sum_{i\in I}$}}
\psfrag{Xalp}{$\overline{\alpha}$} 
\psfrag{Xbeta}{$\overline{\beta}$}
\psfrag{Xgj}{$g_j$} \psfrag{ifinotj}{\large{if $i\not= j$}}
\psfrag{XSnLm}{$S^n_-\otimes \Lambda^m_+$} 
\psfrag{XLmSn}{$\Lambda^m_+\otimes S^n_-$}
\psfrag{XLm1Sn1}{$\Lambda^{m-1}_+\otimes S^{n-1}_-$}
\psfrag{XSumkb}{$-\sum_{b=0}^{k-2} (k-b-1)$}
\psfrag{Xa1}{$\alpha_1$}
\psfrag{Xa2}{$\alpha_2$}
\psfrag{Xb1}{$\beta_1$}
\psfrag{Xb2}{$\beta_2$}
\psfrag{A}{$\mc{A}$}
\psfrag{B}{$\mc{B}$} 
\newpage

%%%%%%%%%%%%%%%%%%%%%%%%%%%%
%%%%%%%%%%%%%%%%%%%%%%%%%%%%
%%
%%   INTRODUCTION 
%%
%%%%%%%%%%%%%%%%%%%%%%%%%%%%
%%%%%%%%%%%%%%%%%%%%%%%%%%%%

\section{Introduction} 

In this paper we propose a graphical calculus for a categorification of 
the Heisenberg algebra. The one-variable Heisenberg algebra has 
generators $p,q$, one defining relation $pq-qp=1$, and appears as
the algebra of operators in the quantization of the  
harmonic oscillator. 
A fundamental role in the quantum field theory is played by its 
infinitely-generated analogue, the algebra with 
generators $p_i,q_i$, for $i$ in some infinite set $I$, and relations 
\begin{equation} \label{eq-heisen} 
p_i q_j = q_j p_i+\delta_{i,j}1, \ \ \ \ 
 p_i p_j=p_j p_i, \ \ \ \ q_i q_j = q_j q_i .
\end{equation}  
In Section~\ref{sec-new} we define a strict monoidal category $\Cone$ with 
two generating objects $Q_+$ and $Q_-$, and morphisms 
between tensor products of these objects given by linear combinations of 
certain planar diagrams modulo local relations. The category is $\Bbbk$-linear 
over a ground commutative ring $\Bbbk$, and we specialize $\Bbbk$ to  
a field of characteristic $0$.  The endomorphism rings of tensor powers 
$Q_+^{\otimes n}$ and $Q_-^{\otimes n}$ contain the group algebra 
$\Bbbk[S_n]$ of the symmetric group. The symmetrization and antisymmetrization 
idempotents in $\Bbbk[S_n]$ produce objects in the Karoubi envelope $\Ctwo$ of 
$\Cone$. These objects can be viewed as symmetric and exterior powers of the 
generating objects  $Q_+$ and $Q_-$. Consequently, we denote them by 
\begin{equation} 
S^n_+:= S^n(Q_+),  \ \  \Lambda^n_+:= \Lambda^n(Q_+), \ \ 
S^n_- :=    S^n(Q_-),  \ \  \Lambda^n_-:=  \Lambda^n(Q_-),
\end{equation} 
and call them the symmetric and exterior powers of $Q_+$ and $Q_-$. When $n=0$, 
$$ S^0_+\cong S^0_-\cong \Lambda^0_+\cong \Lambda^0_-\cong \mathbf{1},$$ 
where $\mathbf{1}$ is the identity object of the monoidal category $\Ctwo$, 
$\mathbf{1}\otimes M = M$ for any $M$. We also set 
$$S^n_+ = S^n_- = \Lambda^n_+ = \Lambda^n_- =0 \ \ \mathrm{if} 
\ \ n<0.$$

\begin{prop} \label{prop-can-intr}
There are canonical isomorphisms in $\Ctwo$
\begin{eqnarray*} 
S^n_- \otimes \Lambda^m_+ & \cong & (\Lambda^m_+ \otimes S^n_-) \oplus 
  (\Lambda^{m-1}_+ \otimes S^{n-1}_-)   , \\ 
S^n_- \otimes S^m_- & \cong & S^m_- \otimes S^n_- , \\ 
\Lambda^n_+ \otimes \Lambda^m_+ & \cong & \Lambda^m_+\otimes \Lambda^n_+. 
\end{eqnarray*}
\end{prop} 

\noindent 
These isomorphisms are constructed in Section~\ref{subsec-Karoubi}. 
Since $\Ctwo$ is monoidal, its Grothendieck group $K_0(\Ctwo)$ is a ring. 
It has generators $[M]$ over objects $M$ of $\Ctwo$ and relations 
$[M_1]=[M_2]+[M_3]$ whenever $M_1\cong M_2\oplus M_3.$ 
The multiplication is defined by $[M_1][M_2] := [M_1\otimes M_2 ]$. 

\begin{cor} \label{cor-rels} The following equalities hold in $K_0(\Ctwo)$: 
\begin{eqnarray*} 
 [ S^n_-] [\Lambda^m_+] & = & [\Lambda^m_+][ S^n_-] +[\Lambda^{m-1}_+]  
  [ S^{n-1}_- ] , \\
  { [ S^n_- ] [ S^m_- ] } & = &  [ S^m_- ]   [ S^n_- ]  ,     \\     
 {[ \Lambda^n_+ ][ \Lambda^m_+ ]} & = & [ \Lambda^m_+ ][ \Lambda^n_+ ]. 
\end{eqnarray*} 
\end{cor} 

Let $H_{\Z}$ be the unital ring with generators 
$a_n, b_n, n\ge 1$ and defining relations
\begin{eqnarray} 
a_n b_m & = & b_m a_n + b_{m-1} a_{n-1},  \label{eq-ab1} \\ 
 a_n a_m & = & a_m a_n, \label{eq-ab2}\\
 b_n b_m & = & b_m b_n. \label{eq-ab3}
\end{eqnarray} 
We simply rewrote relations in Corollary~\ref{cor-rels} using $a_n$ in place of 
 $[S^n_-]$  and $b_m$ instead of $[\Lambda^m_+]$. 
Also set $a_0=b_0=1$, $a_n=b_n=0$ for $n<0$, and require that the above relations 
hold for any $n,m\in\Z$. Any product of $a$'s and $b$'s can be converted into a 
linear combination with nonnegative integer coefficients of monomials 
in $b$'s times monomials in $a$'s,   
\begin{equation}\label{eq-prods}
b_{m_1}b_{m_2} \dots b_{m_k} a_{n_1}a_{n_2} \dots a_{n_r}
\end{equation} 
with $1\le m_1 \le m_2 \le \dots \le m_k,$ $1\le n_1 \le n_2 \le \dots \le n_r$.  
The Bergman diamond lemma~\cite{Berg} tells us that this set of elements is a 
basis of $H_{\Z}$ viewed as a free abelian group.  Let $H=H_{\Z} \otimes \C$ 
be the $\C$-algebra with the same generators and relations as $H_{\Z}$. 

Forming generating functions 
$$ A(t) = 1 + a_1 t +a_2 t^2 + \dots, \ \ \ B(u) = 1+ b_1 u + b_2 u^2 + \dots , $$ 
we can rewrite relations (\ref{eq-ab1}) as 
$$ A(t) B(u) = B(u) A(t) (1 + tu) .$$ 
Let 
$$ \widetilde{A}(t) = 1 + t A'(-t) A(-t) , \ \ \ \widetilde{A}(t) = 1 + \tilde{a}_1 t + 
  \tilde{a}_2 t^2 + \dots   . $$ 
It is easy to check that $\widetilde{a}_1, \widetilde{a}_2, \dots $ generate the 
same subalgebra of $H$ as $a_1, a_2, \dots$, and that  
$$ \widetilde{A}(t) B(u) = B(u) \widetilde{A}(t) + \frac{tu}{1-tu} .$$ 
Coefficients of this equation give us relations (\ref{eq-ab4}) below 
\begin{eqnarray} 
\tilde{a}_n b_m & = & b_m \tilde{a}_n + \delta_{n,m} 1 , \label{eq-ab4} \\ 
 \tilde{a}_n \tilde{a}_m & = & \tilde{a}_m \tilde{a}_n, \label{eq-ab5}\\
 b_n b_m & = & b_m b_n. \label{eq-ab6} 
\end{eqnarray} 
Algebra $H$ is isomorphic to the algebra generated by $\tilde{a}_n$'s, $b_m$'s, 
$n,m > 0$,  
with defining relations (\ref{eq-ab4}) - (\ref{eq-ab6}). 
This allows us to identify $H$ with the Heisenberg algebra  
and $H_{\Z}$ with its integral form.  

\vspace{0.06in} 

Corollary~\ref{cor-rels} gives a ring homomorphism  
\begin{equation}\label{eq-gamma} 
\gamma: H_{\Z} \lra K_0(\Ctwo)
\end{equation}
that takes $a_n$ to $[S^n_-]$ and $b_n$ to $[\Lambda^n_+]$. 

\begin{theorem} \label{thm-inj} Map $\gamma$ is injective.
\end{theorem} 

This theorem is proved in Section~\ref{sub-indres}. 

\begin{conj} \label{conj-iso} Map $\gamma$ is an isomorphism. 
\end{conj} 

If true, this conjecture would allow us to view the additive monoidal 
category $\Ctwo$ as a categorification of the integral form $H_{\Z}$ 
of the Heisenberg algebra. 

\vspace{0.06in} 

The degenerate affine Hecke algebra, which we call degenerate AHA following 
a suggestion of Etingof, was introduced by Drinfeld~\cite{Drin} in the $GL(n)$ 
case and by Lusztig~\cite{Lu} in the general case. Cherednik~\cite{Cher1}
classified finite-dimensional irreducible representations of degenerate AHA; its 
centralizing properties were studied by him and Olshanski in~\cite{Cher2,OL}. 
We denote by $DH_n$ the degenerate AHA for $GL(n)$, over the base field $\Bbbk$. 
Under the canonical homomorphism~\cite{Cher1, Drin} from $DH_n$ 
to the group algebra $\Bbbk[S_n]$ of the symmetric group polynomial 
generators of $DH_n$ go to the Jucys-Murphy elements. Okounkov and 
Vershik~\cite{OV1},~\cite{OV2} presented a detailed 
derivation of the basic representation theory of the symmetric group via 
these elements, see also~\cite[Chapter 2]{Kbook},~\cite{CST,DG}. For some other 
uses of Jucys-Murphy's elements and degenerate AHA we refer the reader 
to~\cite{Hora, MO, Ok1, Sn}. 

\vspace{0.06in} 

We will prove in Section~\ref{sec-sizeofC} that the ring of endomorphisms of the object 
$Q_+^{\otimes n}$ in our category is  isomorphic to the tensor product of 
$DH_n$ and the polynomial algebra in infinitely many variables. 
Thus, the degenerate AHA for $GL(n)$ emerges naturally in our approach as 
part of a larger structure. Polynomial generators of $DH_n$ 
acquire graphical interpretation in our calculus as right-twisted curls on strands. 
We also describe a basis, given diagrammatically, of vector spaces of 
morphisms between arbitrary 
tensor products of generators $Q_+$ and $Q_-$ of $\Cone$.

To prove our results, we construct a family of functors from $\Cone$ to the 
category $\Sone$ whose objects are compositions of induction and restriction 
functors between group algebras $\Bbbk[S_n]$ of the symmetric group, and 
morphisms are natural transformations between these functors. The image 
under these functors  of the endomorphism of $Q_+$ given by the 
right curl diagram is the Jucys-Murphy element. The image of the  
counterclockwise circle diagram with $k$ right curls is the $k$-th moment 
of the Jucys-Murphy element. Products of these moments were investigated 
in~\cite{Hora, Ok1, Rt, Sn} in relation to the asymptotic representation 
theory of the symmetric group and free probability. 
It also appears that our construction should be related to the circle of ideas 
considered by Guionnet, Jones and Shlyakhtenko~\cite{GJS} that intertwine  
planar algebras and free probability. 
In addition, one would hope for a relation between our calculus 
and the geometrization of the Heisenberg algebra via Hilbert schemes by 
Nakajima \cite{Nak1}, \cite{Nak2} and Grojnowski~\cite{Groj}, and for a  
possible link with Frenkel, Jing and Wang~\cite{FJW}. 

\vspace{0.06in} 

We discovered monoidal category $\Cone$ by considering compositions of induction and 
restriction functors for standard inclusions of symmetric group algebras 
$\Bbbk[S_n] \subset \Bbbk[S_{n+1}].$ Induction and restriction 
functors for inclusions of finite groups are biadjoint, and 
biadjointness natural transformations can be 
depicted via cap and cup planar diagrams. Furthermore, the composition 
of two induction functors admits a natural endotransformation given by the 
right multiplication by the transposition $(n+1,n+2)$, an endomorphism 
of $\Bbbk[S_{n+2}]$ viewed as 
left $\Bbbk[S_{n+2}]$-module and right $\Bbbk[S_n]$-module. 
We denote this natural transformation by the crossing of two upward-oriented 
strands. Relations on compositions of the crossing, cup, and cap transformations 
that hold for all $n$ (universal relations) are given by equations 
(\ref{eq-mainrels1})-(\ref{eq-mainrels3}). These relations together with 
isotopies of diagrams are exactly the 
defining relations for the additive monoidal category $\Cone$. 

\vspace{0.06in} 

Cautis and Licata~\cite{CL} introduced graded relatives of $\Cone$ and $\Ctwo$ 
associated to finite subgroups of $SU(2)$, identified their Grothendieck rings with 
certain quantized Heisenberg algebras, and constructed an action of their categories 
on derived categories of coherent sheaves on Hilbert schemes of points on the ALE spaces. 
Hom spaces in Cautis-Licata monoidal categories 
carry a natural grading (absent in our case) with finite-dimensional 
homogeneous terms and vanishing negative degree homs on certain objects, 
leading to a proof that their analogue of the map $\gamma$ is an isomorphism.       

\vspace{0.06in} 

By themselves, Heisenberg algebras are rather simple constructs. 
Their value is in the structures that  quantum field 
theory builts on top of them, for instance, the structures of vertex operator algebras. 
The problem posed by Igor Frenkel~\cite{Fre} to categorify just the simplest 
vertex operator algebra remains wide open - perhaps our paper will serve as  
a small step towards this goal. 

\vspace{0.06in} 

{\bf Acknowledgments:} On a number of occasions Igor Frenkel emphasized to the 
author the importance of categorifying various structures related to the 
symmetric functions and Heisenberg algebras. 
More recently, Tony Licata asked directly if there exists a categorification 
of the Heisenberg algebra, and this paper was inspired by his question. 
The author would like to acknowledge partial support by the National Science Foundation 
via grants DMS-0706924 and DMS-0739392.

%%%%%%%%%%%%%%%%%%%%%%%%
%%%%%%%%%%%%%%%%%%%%%%%%
%%
%%  CALCULUS
%%
%%%%%%%%%%%%%%%%%%%%%%%%
%%%%%%%%%%%%%%%%%%%%%%%%

\section{A new graphical calculus} \label{sec-new}

\subsection{Local moves and their consequences} \label{subsection-moves}

Fix a commutative ring $\Bbbk$ and consider the following additive $\Bbbk$-linear 
monoidal category $\Cone$ generated by two  objects $Q_+$ and $Q_-$. An object 
of $\Cone$ is a finite direct sum of  
tensor products $Q_{\epsilon_1}\otimes \dots \otimes Q_{\epsilon_m}$, 
denoted $Q_{\epsilon},$ where $\epsilon=\epsilon_1\dots \epsilon_m$ are 
finite sequences of signs. Thus, $Q_{\epsilon\epsilon'}\cong Q_{\epsilon}\otimes 
Q_{\epsilon'}$ for sequences $\epsilon$, $\epsilon'$ and their concatenation 
$\epsilon\epsilon'$. The unit object corresponds to the  
empty sequence: $\mathbf{1}=Q_{\emptyset}.$ 

The space of homomorphisms $\Hom_{\Cone}(Q_{\epsilon},Q_{\epsilon'})$
for sequences $\epsilon,  \epsilon'$ is the $\Bbbk$-module generated by 
suitable planar diagrams, modulo local relations. The diagrams consist of oriented 
compact one-manifolds immersed into the plane strip $\R\times [0,1]$, modulo 
rel boundary isotopies. The relations are 

\begin{equation}\label{eq-mainrels1}
{\psfig{figure=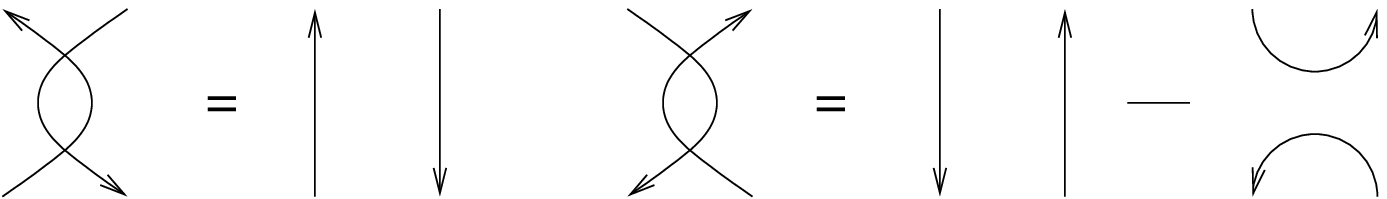,height=2cm}}\end{equation}

\begin{equation}\label{eq-mainrels2}
{\psfig{figure=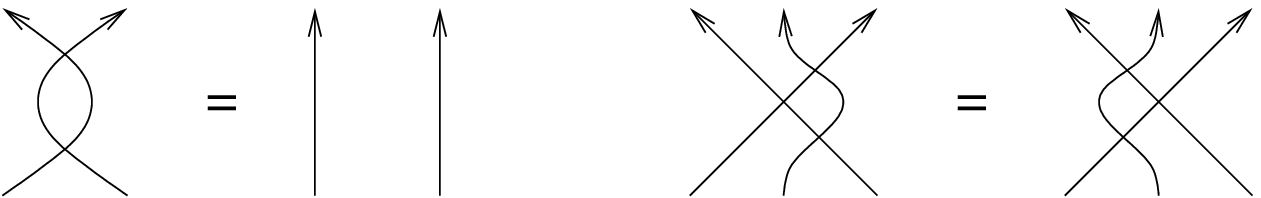,height=2cm}}\end{equation}

\begin{equation}\label{eq-mainrels3}
{\psfig{figure=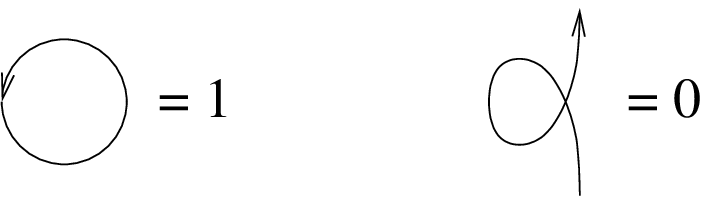,height=2cm}}\end{equation}

We require that the endpoints of the one-manifold are located at 
$\{1,\dots, m\}\times \{0\}$ and $\{1,\dots, k\}\times \{1\}$ and call 
these the lower and upper endpoints, respectively, where $m$ and $k$ are 
lengths of sequences $\epsilon$ and $\epsilon'$. Moreover, orientation of the 
one-manifold at the endpoints must match the signs in the sequences 
$\epsilon$ and $\epsilon'$. For instance, the diagram

\begin{center}{\psfig{figure=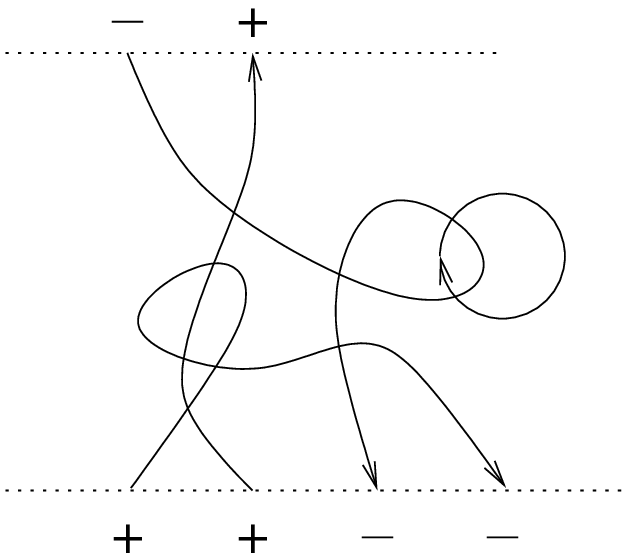,height=5cm}}\end{center}
 
\noindent 
is a morphism from $Q_{++--}$ to $Q_{-+}$. A diagram without 
endpoints is an endomorphism of $\mathbf{1}$. Composition of morphisms 
is given by concatenating the diagrams. The sequence of $n$ pluses is 
denoted $+^n$, the sequence of $n$ minuses $-^n$. 

We have the Heisenberg relation 
$$ Q_{-+}\cong Q_{+-} \oplus \mathbf{1} .$$ 
This isomorphism is canonical and comes from the maps  
between these objects encoded in the following diagram 

\begin{equation} \label{eq-squarepic} 
{\psfig{figure=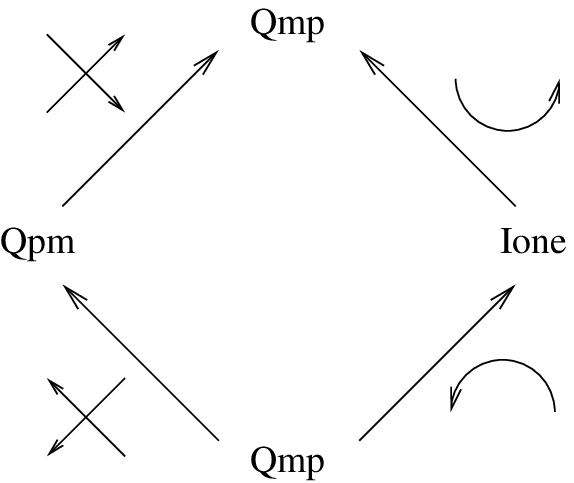}}
\end{equation} 

The four arrows are given by four morphisms, two of which are 
crossings and two are U-turns. The condition that these maps describe a 
decomposition of $Q_{-+}$ as the direct sum of objects $Q_{-+}$ 
and $\mathbf{1}$ is equivalent to relations
(\ref{eq-mainrels1}) and (\ref{eq-mainrels3}), 
modulo the condition that an isotopy of a diagram does 
not change the morphism. The latter condition is 
equivalent to the biadjointness of functors of tensoring with $Q_+$ and 
$Q_-$, with the biadjointness transformations given by the four $U$-turns 

\vspace{0.1in} 
\begin{center}{\psfig{figure=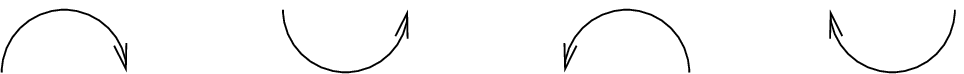}}\end{center}
\vspace{0.1in} 

\noindent 
(see Section~\ref{subsec-biadjoint} for details).  

Moving the lower endpoints of a diagram up via a multiple cups diagram 
leads to canonical isomorphisms 
\begin{equation} \label{eq-cis}
 \Hom_{\Cone}(Q_{\epsilon}, Q_{\epsilon'}) \cong 
   \Hom_{\Cone}({\mathbf 1}, Q_{\overline{\epsilon}\epsilon'}) \cong 
   \Hom_{\Cone}({\mathbf 1}, Q_{\epsilon'\overline{\epsilon}}), 
\end{equation}  
related to the biadjointness of tensoring with  $Q_{\epsilon}$ and  
$Q_{\overline{\epsilon}}$. 
Here $\overline{\epsilon}$ is the sequence $\epsilon$ with 
the order and all signs reversed. 
Biadjointness natural transformations satisfy the cyclicity 
condition~\cite{Bartlett1, Bartlett2,CKS, AL1}, 
which follows at once from the definition of $\Cone$. 

\vspace{0.1in} 

The two relations in (\ref{eq-mainrels1}) allow simplification 
of a double crossing for oppositely oriented intervals. 
The first relation in (\ref{eq-mainrels3}) says that a counterclockwise oriented 
circle equals one. 
Thus, an innermost counterclockwise circle can be erased from the diagram without 
changing the value of the diagram viewed 
as an element of the hom space between the functors. 

There are two possible types of curls on strands: a left  curl 
$\raisebox{-0.6cm}{\psfig{figure=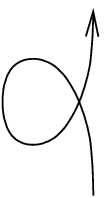,height=1.5cm}}$ and 
a right  curl $\raisebox{-0.6cm}{\psfig{figure=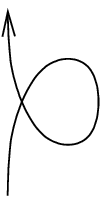,height=1.5cm}}.$
The second relation in (\ref{eq-mainrels3})
says that a diagram that contains a left  curl subdiagram is zero. 

Defining local relations in $\Cone$ imply the following relations 

\begin{center}
{\psfig{figure=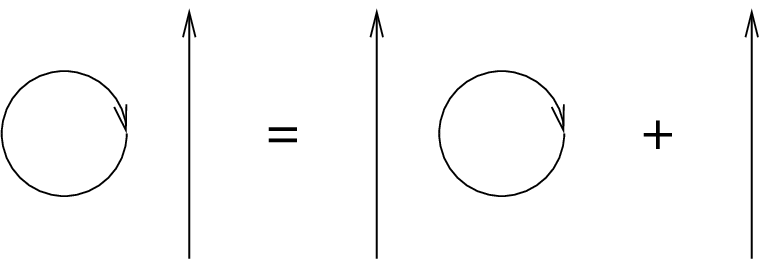,height=1.6cm}}, \hspace{0.4in} 
{\psfig{figure=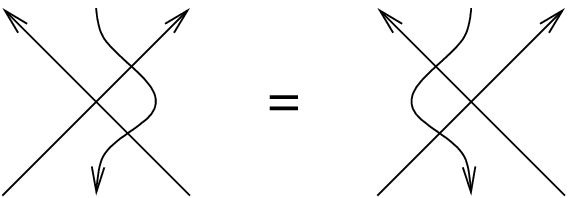,height=1.6cm}}.
\end{center}

The second relation, jointly with the original ones, implies that 
the triple intersection move holds for any orientation of the 3 strands: 

\begin{center}
{\psfig{figure=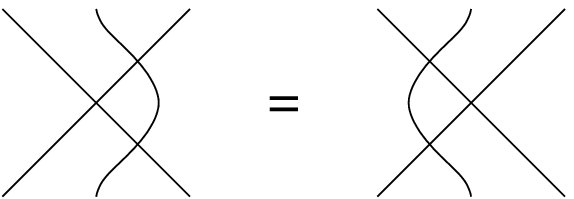,height=1.6cm}} \end{center}

Furthermore, right curls 
can be moved across intersection points, modulo simpler diagrams: 

\begin{center} 
{\psfig{figure=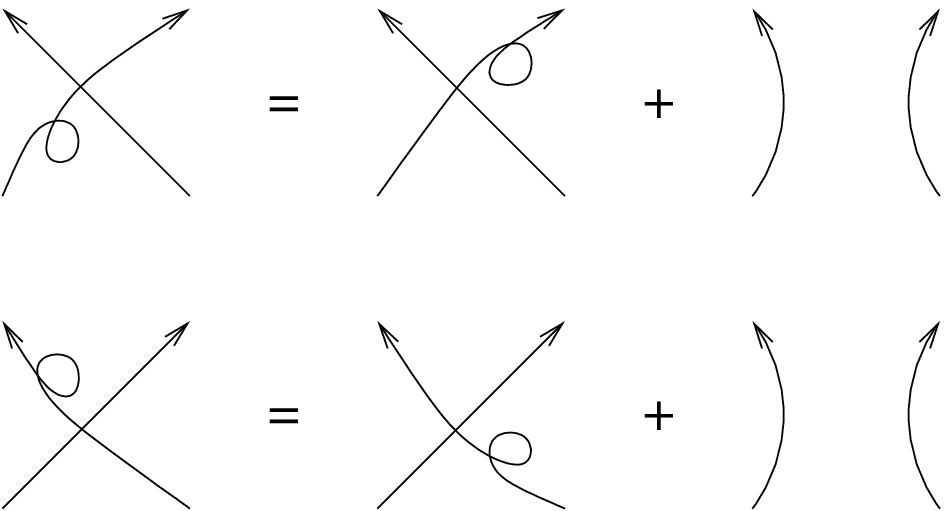,height=3.5cm}}
\end{center}

It will be convenient to denote a right curl by a dot on a strand, 
and $k$-th power of a right curl by a dot with $k$ next to it: 

\begin{center}
{\psfig{figure=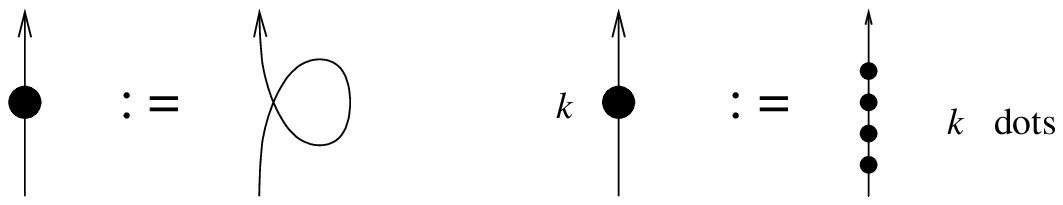}}
\end{center}

The above relations can be rewritten as 

\begin{center} 
{\psfig{figure=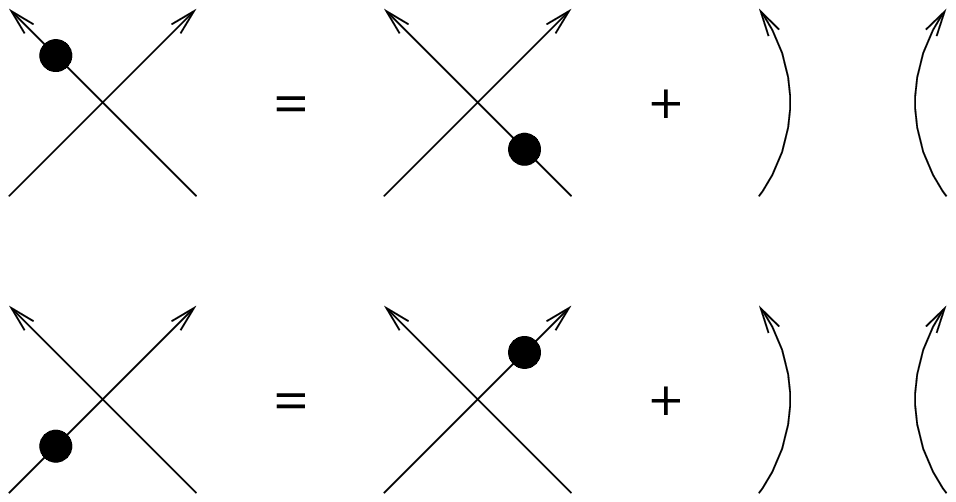,height=3cm}}
\end{center}

Together with the earlier ones, these relations show that there 
is a homomorphism from the degenerate affine Hecke algebra 
$\mbox{DH}_n$ with coefficients in $\Bbbk$ to the $\Bbbk$-algebra of 
endomorphisms of the object $Q_{+^n}$. The permutation 
generator $s_i$ of $\mbox{DH}_n$ goes to the permutation diagram 
of the $i$-th and $i+1$-st strands, and the polynomial generator $x_i$ goes 
to the dot on the $i$-th strand: 

\begin{center} 
{\psfig{figure=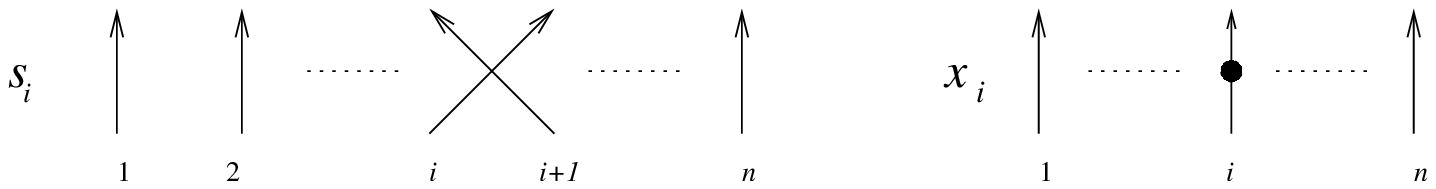,height=1.6cm}}
\end{center}

\vspace{0.1in} 

Note that $\raisebox{-0.3cm}{\psfig{figure=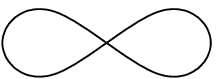}}=0.$
Indeed, the figure eight diagram, for any orientation, contains both 
left and right curls, and, therefore, equals to zero in our calculus. 

A strand with $k$ dots can be closed
into either a clockwise-oriented or a  counterclockwise-oriented circle with $k$ dots. 
Denote these circles by $c_k$ and $\tilde{c}_k$, respectively: 

\begin{equation*} c_k := \ \  \raisebox{-0.3cm}{\psfig{figure=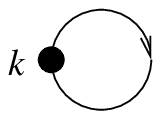}}\ , 
 \ \ \ \ \  \tilde{c}_k :=  \ \ \raisebox{-0.3cm}{\psfig{figure=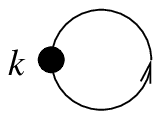}}\ . 
\end{equation*} 

Counterclockwise circles can be expressed as linear combinations of products of 
clockwise circles. For the first few values of $k$, these are 

\begin{eqnarray*}
  \tilde{c}_0 & = & 1, \\
  \tilde{c}_1 & = & 0, \\
  \tilde{c}_2 & = & c_0, \\
  \tilde{c}_3 & = & c_1  ,\\
  \tilde{c}_4 & = & c_2 + c_0^2 , \\
  \tilde{c}_5 & = & c_3 + 2 c_1 c_0.  
\end{eqnarray*} 

These equations are obtained by expanding each dot into a left curl and 
then operating on the resulting diagram via the rules of the graphical calculus. 
A counterclockwise circle with one dot expands into the figure eight diagram,  which 
is $0$. For another example, 

\begin{equation*} 
\tilde{c}_3 = \ \ \raisebox{-0.7cm}{\psfig{figure=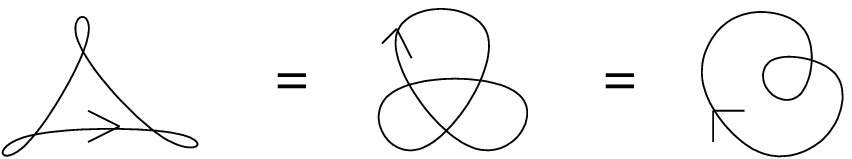}} \ \ \ = 
 c_1. 
\end{equation*} 

\begin{prop} For $k>0$ we have  
\begin{equation} \label{eq-tildec}  
 \tilde{c}_{k+1} = \sum_{a=0}^{k-1} \tilde{c}_a c_{k-1-a}. 
\end{equation} 
\end{prop}  

\emph{Proof} is  the following computation:

\begin{center}
{\psfig{figure=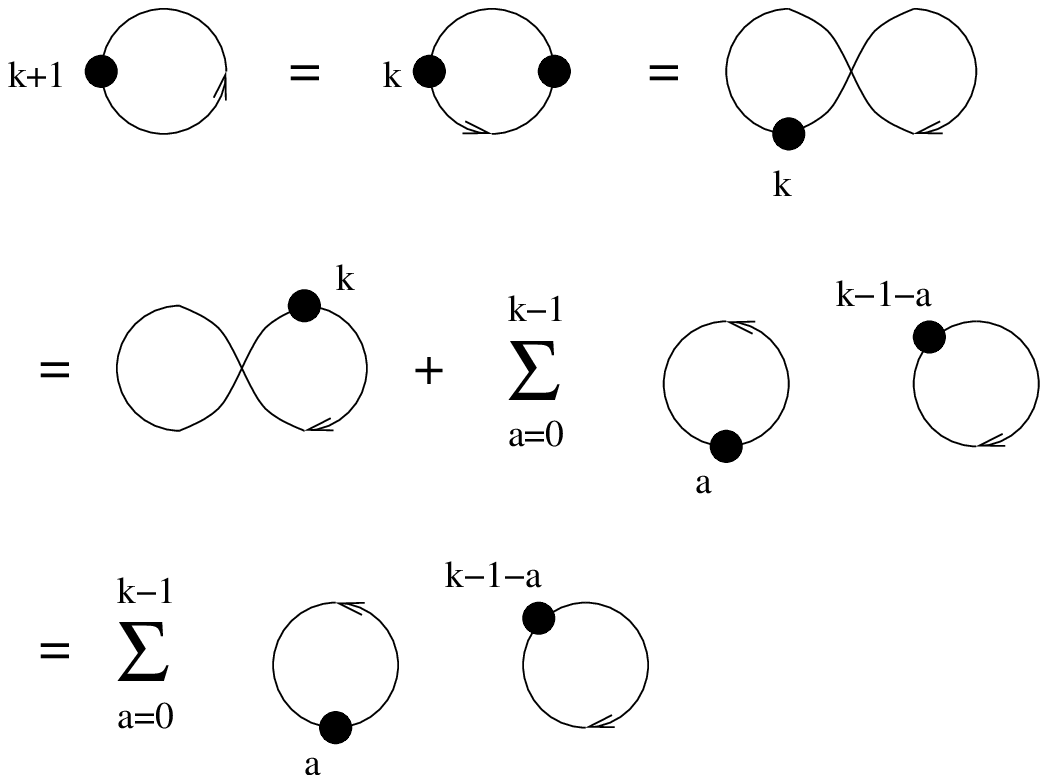,height=7cm}}
\end{center}
In the second equality, we converted a dot to a right curl, and, in the 
third equality, moved $k$ dots through a crossing. The first term 
on the second line equals $0$, since it contains a left curl.
$\square$ 

\vspace{0.06in} 

Iterating this formula, one obtains an expession for $\tilde{c}_k$ as 
a polynomial function of $c_m$, $m\le k-2$. 
Vice versa, each $c_m$ can be written as a polynomial in $\tilde{c}_k, $ 
$k\le m+2$. Let $t$ be a formal variable and write 

$$c(t) = \sum_{i=0}^{\infty} c_i t^i, \ \ \ \ 
   \tilde{c}(t) =  \sum_{i=0}^{\infty} \tilde{c}_i t^i.$$ 
Formula (\ref{eq-tildec}) turns into $ t^2 c(t) \tilde{c}(t) = \tilde{c}(t) - 1,$
so that 
$$ \tilde{c}(t) = \frac{1}{1 - t^2 c(t)}.$$ 

\vspace{0.06in} 

The following identities, called \emph{bubble moves} by analogy with~\cite{AL1}, 
hold.

\begin{center} 
{\psfig{figure=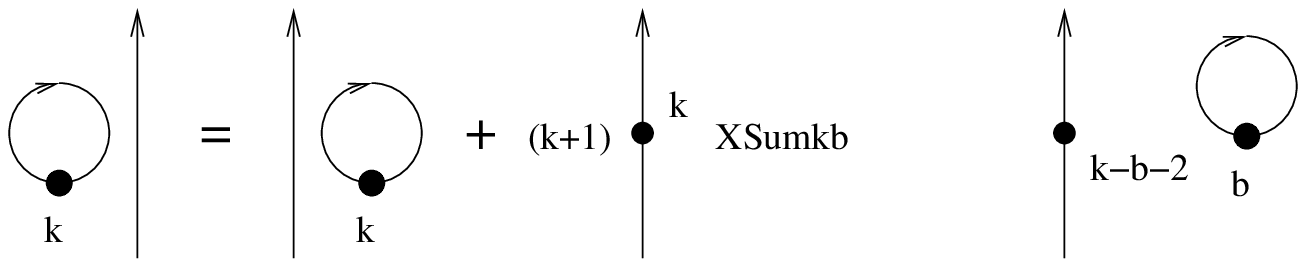}}
\end{center}

\begin{center} 
{\psfig{figure=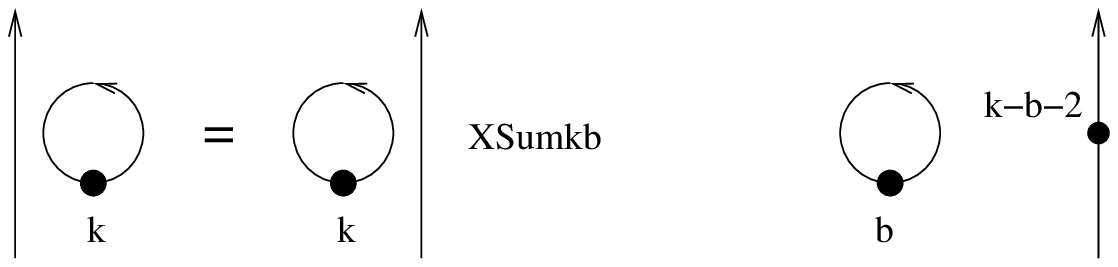}}
\end{center}

A closed diagram $D$ defines an endomorphism of the object $\mathbf{1}$ 
of $\Cone$. Using the local moves,  
such diagram $D$ can be converted into a  linear combinations of crossingless 
diagrams that consist of nested dotted circles. Furthermore, bubble moves 
can be used to split apart nested circles. Lastly, convert counterclockwise circles 
into linear combinations of products of clockwise circles. Therefore, 
a closed diagram can be written as a linear combination of products of 
dotted clockwise circles. We see 
that the endomorphism algebra $\End_{\Cone}(\mathbf{1})$ is 
a quotient of the polynomial algebra $\Pi:=\Bbbk[c_0,c_1,c_2, \dots ]$ 
in countably many variables via the map 
\begin{equation} \label{eq-map-iso}
\psi_0  \ : \ \Pi=\Bbbk[c_0,c_1,c_2, \dots ] \lra \End_{\Cone}(\mathbf{1})
\end{equation}
that takes $c_k$ to the clockwise circle with $k$ dots (we allowed 
ourselves the liberty of using $c_k$ to denote both a formal variable 
and its image in the endomorphism algebra). 

\begin{prop} \label{prop-psi-iso} Map $\psi_0$ is an isomorphism. 
\end{prop} 

This proposition will be proved in Section~\ref{sec-sizeofC}. 

The endomorphism algebra of $Q_{+^m}$ is spanned 
by all diagrams that have $m$ upper and $m$ lower endpoints and such that 
at each endpoint the strand is oriented upward. 
A homomorphism from the degenerate affine Hecke algebra $\mbox{DH}_m$ 
to $\End_{\Cone}(Q_{+^m})$ was described earlier. Placing 
a closed diagram to the right of a diagram representing an element of $\mbox{DH}_m$ 
gives a homomorphism 
\begin{equation}\label{eq-psim} 
 \psi_m \ : \\  \mbox{DH}_m \otimes \Pi \lra \End_{\Cone}(Q_{+^m}).
\end{equation}  
It is easy to see that $\psi_m$ is surjective, by taking a diagram representing 
an element on the right hand side, and inductively simplifying it to a linear combination 
of diagrams that come from a standard basis of the left hand side. 
Namely, any diagram representing an endomorphism of $Q_{+^m}$ 
is a combination of diagrams that consist of a  permutation $\sigma \in S_m$, 
some number (possibly zero) of dots on each strand above the permutation
diagram,  and a monomial in dotted clockwise circles to the right: 

\begin{center} 
{\psfig{figure=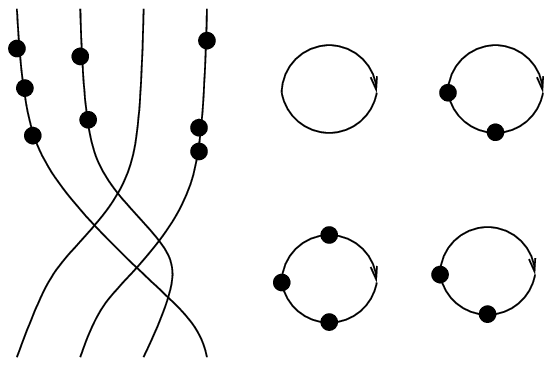}}
\end{center}

We write these basis elements as 
$x_1^{a_1}\dots x_m^{a_m}\cdot \sigma \cdot c_0^{b_0} c_1^{b_1}\dots c_k^{b_k}.$
In the above example, the element is 
$x_1^3x_2^2x_4^3 \cdot (1324)\cdot c_0 c_2^2 c_3.$ 

Surjectivity of $\psi_m$ can be strengthened to the following result. 

\begin{prop} \label{prop-psim-iso}
$\psi_m$ is an isomorphism. \end{prop} 

Injectivity of $\psi_m$  is proved in Section~\ref{sec-sizeofC}.

\vspace{0.06in} 

We now describe a spanning set in the $\Bbbk$-module 
$\Hom_{\Cone}(Q_{\epsilon}, Q_{\epsilon'})$ for any 
sequences $\epsilon,\epsilon'$. Let $k$ be the total number of 
pluses in these two sequences. This hom space is nontrivial only 
if the total number of minuses in these two sequences is $k$ as 
well, which we'll assume from now on to be the case. 
The spanning set, denoted $B(\epsilon, \epsilon')$, is obtained 
by forming all possible oriented matchings of these two points 
via $k$ oriented segments in the plane strip $\R\times [0,1]$. 
We assume that sequences $\epsilon$ and $\epsilon'$ are 
written at the bottom and top of the strip, the segments are 
embedded in the strip, and their orientations at the endpoints 
match corresponding elements of $\epsilon$ and $\epsilon'$. 
Each two segments intersect at most once, and no triple intersections 
are allowed, see an example below for $\epsilon=--+$ and $\epsilon'=+---+$.

\begin{center}{\psfig{figure=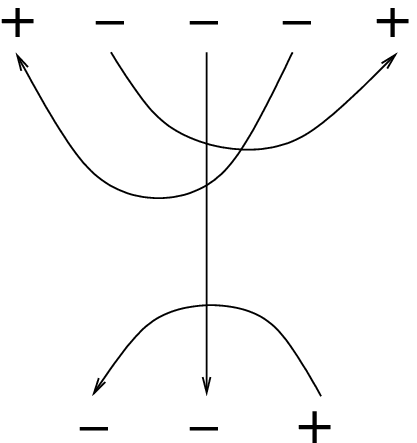}}\end{center}

Select an interval disjoint from intersections near the out endpoint 
of each interval and put any number (perhaps 0) of dots on it. In the rightmost 
region of the diagram, draw some number of clockwise oriented 
disjoint nonnested circles with no dots, some number of such circles with 
one dot, two dots, etc., with finitely-many circles  in total. 
The resulting set of diagrams $B(\epsilon, \epsilon')$ is parametrized 
by $k!$ possible matchings of the $2k$ oriented endpoints, by a sequence of 
$k$ nonnegative integers describing the number of dots on each 
interval, and by a finite sequence of nonnegative integers listing 
the number of clockwise oriented circles with no dots, one dot, and so on. 
An example of a diagram in $B(--+, +---+)$ is depicted below.

\begin{center}{\psfig{figure=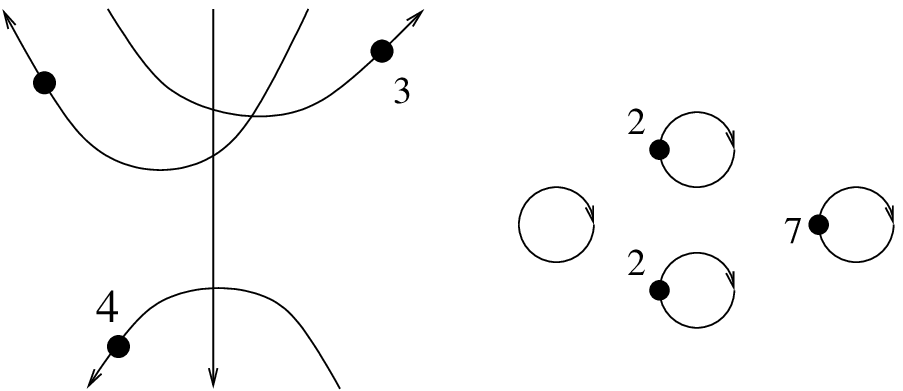}}\end{center}

In this picture, we put no dots, one dot, three dots, and four dots 
on the four arcs of the matching, added one bubble with no dots, 
two with two dots and one with five dots.   

It is rather straightforward to check that $B(\epsilon,\epsilon')$ is 
a spanning set of the $\Bbbk$-vector space 
$\Hom_{\Cone}(Q_{\epsilon}, Q_{\epsilon'})$. 

\begin{prop} \label{prop-psi-iso2}
For any sign sequences $\epsilon,\epsilon'$ the set $B(\epsilon, \epsilon')$ 
constitutes a basis of the $\Bbbk$-vector space 
$\Hom_{\Cone}(Q_{\epsilon}, Q_{\epsilon'}).$ 
\end{prop} 

Thus, we claim that the set $B(\epsilon, \epsilon')$ is also linearly independent. 
This proposition holds as well for $\Bbbk$ being any commutative ring rather 
than a field, with $B(\epsilon, \epsilon')$ being a basis of the free $\Bbbk$-module 
$\Hom_{\Cone}(Q_{\epsilon}, Q_{\epsilon'}).$

Notice that Proposition~\ref{prop-psim-iso} is a special case of this proposition, for 
$\epsilon=\epsilon'=+^m$. Proposition~\ref{prop-psi-iso2} follows 
from proposition~\ref{prop-psim-iso}, functor isomorphisms 
$Q_{-+}\cong Q_{+-}\oplus \Id$ and arguments 
similar to the ones in~\cite[Section 2.2]{KL3}. First, canonical isomorphisms 
(\ref{eq-cis}) take the set $B(\epsilon,\epsilon')$ to $B(\emptyset, 
\overline{\epsilon}\epsilon')$ and $B(\emptyset, \epsilon'\overline{\epsilon})$, 
respectively, and it is then enough to show that $B(\emptyset, \epsilon)$ 
is linearly independent for any sequence $\epsilon$ with $k$ pluses and 
$k$ minuses. 

Proposition~\ref{prop-psim-iso} implies linear independence 
for $k=0, 1$, and for the sequence $+^k -^k$ for any $k$. 
Assume that $\epsilon = \epsilon_1 - + \ \epsilon_2$ for 
some sequences $\epsilon_1, \epsilon_2$. Assume by induction on 
$k$ and by induction on the lexicographic order among length $2k$ 
sequences that the sets $B(\emptyset, \epsilon_1\epsilon_2)$ 
and $B(\emptyset, \epsilon_1+- \ \epsilon_2)$ are linearly independent 
in their respective hom spaces. 

Two upper arrows in 
the diagram (\ref{eq-squarepic}) lead to a canonical decomposition 
$$ Q_{\epsilon_1 - + \epsilon_2} \cong Q_{\epsilon_1 + - \epsilon_2} 
\oplus Q_{\epsilon_1 \epsilon_2}.$$ 
Under this isomorphism sets $B(\emptyset, \epsilon_1\epsilon_2)$ 
and $B(\emptyset, \epsilon_1+- \ \epsilon_2)$ get mapped to two 
subsets of $\Hom_{\Cone} ({\mathbf 1}, Q_{\epsilon_1 - + \epsilon_2} )$. 
Denote by $B$ the union of these two subsets. It is easy to see that 
linear independence of $B(\emptyset,\epsilon_1 - + \ \epsilon_2)$ 
is equivalent to linear independence of $B$, which we know by induction. 
Proposition~\ref{prop-psi-iso2} follows. $\square$ 

\vspace{0.2in}

By the thickness of a diagram in $B(\epsilon, \epsilon')$ we call 
the number of arcs connecting lower and upper endpoints. 
The diagram depicted earlier has thickness one. 
For each $k$ and $\epsilon$, the subset of diagrams of thickness at most $k$ 
is a 2-sided ideal in the endomorphism ring of $Q_{\epsilon}$, 
since thickness cannot increase upon composition.  
For $\epsilon= +^n -^m$ and $k=n+m-1$ we denote the 
corresponding ideal by $J_{n,m}$. It is spanned by diagrams with 
at least one arc connecting a pair of upper endpoints (and, necessarily, 
at least one arc connecting a pair of lower endpoints). 
It is easy to see that the quotient of the endomorphism 
ring of $Q_{+^n -^m}$ by this ideal is naturally isomorphic 
to the tensor product $DH_n\otimes DH_m^{op} \otimes \Pi$, and 
the short exact sequence 
\begin{equation} \label{eq-splits} 
0 \lra J_{n,m} \lra \End_{\Cone}(Q_{+^n -^m}) \lra DH_n \otimes 
DH_m^{op} \otimes \Pi \lra 0 
\end{equation} 
admits a canonical splitting. Notice also that $DH_m^{op}\cong DH_m$. 

\vspace{0.2in} 

We now list some obvious symmetries of $\Cone$. 
The map that assigns  $(-1)^{w(D)} D$ to a diagram $D$, where $w(D)$ is 
the number of crossings plus the number of dots of $D$, 
extends to an involutive autoequivalence $\xi_1$ of $\Cone$. We have $\xi_1^2=\Id$ 
(equality and not just isomorphism). Autoequivalence $\xi_1$ exchanges 
$S^n_+$ with $\Lambda^n_+$ and $S^n_-$ with $\Lambda^n_-$. 

\vspace{0.06in} 

Denote by $\xi_2$ the symmetry of category $\Cone$ given on diagrams
by  reflecting about the x-axis and reversing orientation. This symmetry is an 
involutive monoidal contravariant autoequivalence of $\Cone$. 

\vspace{0.06in} 

Denote by $\xi_3$ the symmetry of category $\Cone$ given on diagrams 
by reflecting about the y-axis and reversing orientation. This symmetry is an 
involutive antimonoidal autoequivalence of $\Cone$.  Being antimonoidal means 
that it reverses the order of elements in the tensor product:  
$\xi_3(M\otimes N) = \xi_3(N) \otimes \xi_3(M).$ 

\vspace{0.06in} 

Symmetries $\xi_1, \xi_2, \xi_3$ pairwise commute and generate an action 
of $(\Z/2)^3$.

%%%%%%%%%%%%%%%%%%%%%%%%%
%%%%%%%%%%%%%%%%%%%%%%%%
%%
%%  KAROUBI ENVELOPE
%%
%%%%%%%%%%%%%%%%%%%%%%%%
%%%%%%%%%%%%%%%%%%%%%%%%

\subsection{Karoubi envelope and projectors}\label{subsec-Karoubi} 

The two relations in (\ref{eq-mainrels2}) tell us that 
upward oriented crossings 
satisfy the symmetric group relations and give us a canonical homomorphism 
$$ \Bbbk [S_n] \lra \End_{\Cone}(Q_{+^n}) $$ 
from the group algebra of the symmetric group to the endomorphism 
ring of the $n$-th tensor power of $Q_+$. Turning the diagrams by $180$
degrees, we obtain a canonical homomorphism 
$$ \Bbbk [S_n] \lra  \End_{\Cone}(Q_{-^n}).$$ 
Assume that $\Bbbk$ is a field of characteristic $0$. Then we can use symmetrizers 
and antisymmetrizers, and, more generally, Young symmetrizers, to 
produce idempotents in $\End_{\Cone}(Q_{+^n})$. At this 
point it is convenient to introduce the Karoubi envelope of $\Cone,$ 
the category $\Ctwo$ whose objects 
are pairs $(P,e)$, where $P$ is an object of $\Cone$ and $e: P\lra P$ is 
an idempotent endomorphism, $e^2=e$. Morphisms from $(P,e)$ to $(P',e')$ are 
maps $f: P \lra P'$ in $\Cone$ such that $e'fe=f$. It is immediate 
that $\Ctwo$ is a $\Bbbk$-linear additive monoidal category. 

To the complete symmetrizer 
$$e(n)\in \Bbbk[S_n], \ \ \  e(n)=\frac{1}{n!} \sum_{\sigma\in S_n} \sigma$$  
we assign the object $S^n_+:=(Q_{+^n},e(n))$ in $\Ctwo$. 
Following Cvitanovi\'c~\cite{CVBook}, which contains diagrammatics 
for Young symmetrizers and antisymmetrizers,  
we depict $S^n_+$ as a white box labelled $n$. The inclusion morphism 
$S^n_+ \lra Q_{+^n}$ is depicted by a white box with $n$ upward oriented 
lines emanating from the top. The projection $Q_{+^n}\lra S^n_+$
is depicted by a white box with $n$ upward oriented lines at the bottom. 
The composition $Q_{+^n}\lra S^n_+\lra Q_{+^n}$ is depicted likewise. 

\begin{center} {\psfig{figure=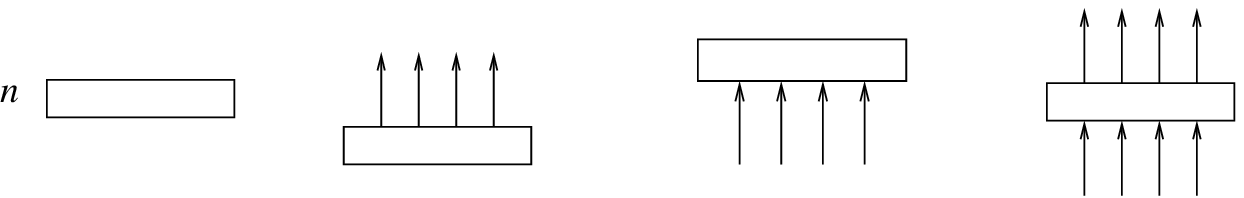}}\end{center}

To the complete antisymmetrizer 
$$e'(n)\in \Bbbk[S_n], \ \ \  e'(n)=\frac{1}{n!} \sum_{\sigma\in S_n}
\mathrm{sign}(\sigma) \sigma$$ 
we assign the object $\Lambda^n_+:=(Q_{+^n},e'(n))$ in $\Ctwo$
and depict it and related inclusions and projections to and from $Q_{+^n}$ 
by black boxes with up arrows 

\begin{center} {\psfig{figure=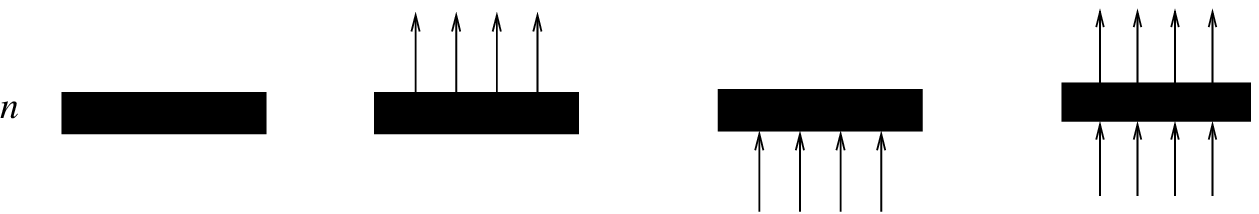}}\end{center}

Define the objects $S^-_n:=  (Q_{-^n},e(n))$  and  
$\Lambda^-_n:= (Q_{-^n},e'(n)) $ as the subobjects of $Q_{-^n}$ 
associated to the symmetrizer $e(n)$ and the antisymmetrizer $e'(n)$ idempotents, 
respectively, under the canonical homomorphism 
$\Bbbk [S_n] \lra  \End_{\Cone}(Q_{-^n})$. We draw $S^-_n$ and 
$\Lambda^-_n$ as white, respectively black, 
boxes, but with the lines at the boxes oriented downward.

\vspace{0.1in} 

We plan to develop the graphical calculus of these diagrams elsewhere. 
Part of the calculus that deals with the lines oriented only upwards (or only 
downwards) is  the graphical calculus of symmetrizers and antisymmetrizers in 
the symmetric group, and can be found in~\cite{CVBook}. 
The latter calculus implies the second and third isomorphisms 
from Proposition~\ref{prop-can-intr} in the introduction. 
For instance, the second isomorphism is realized by the diagram 

\begin{center} {\psfig{figure=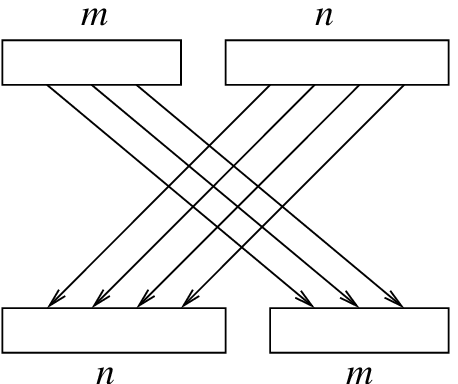}}\end{center}

The first family of isomorphisms in Proposition~\ref{prop-can-intr} 
is realized by the maps 

\begin{center} {\psfig{figure=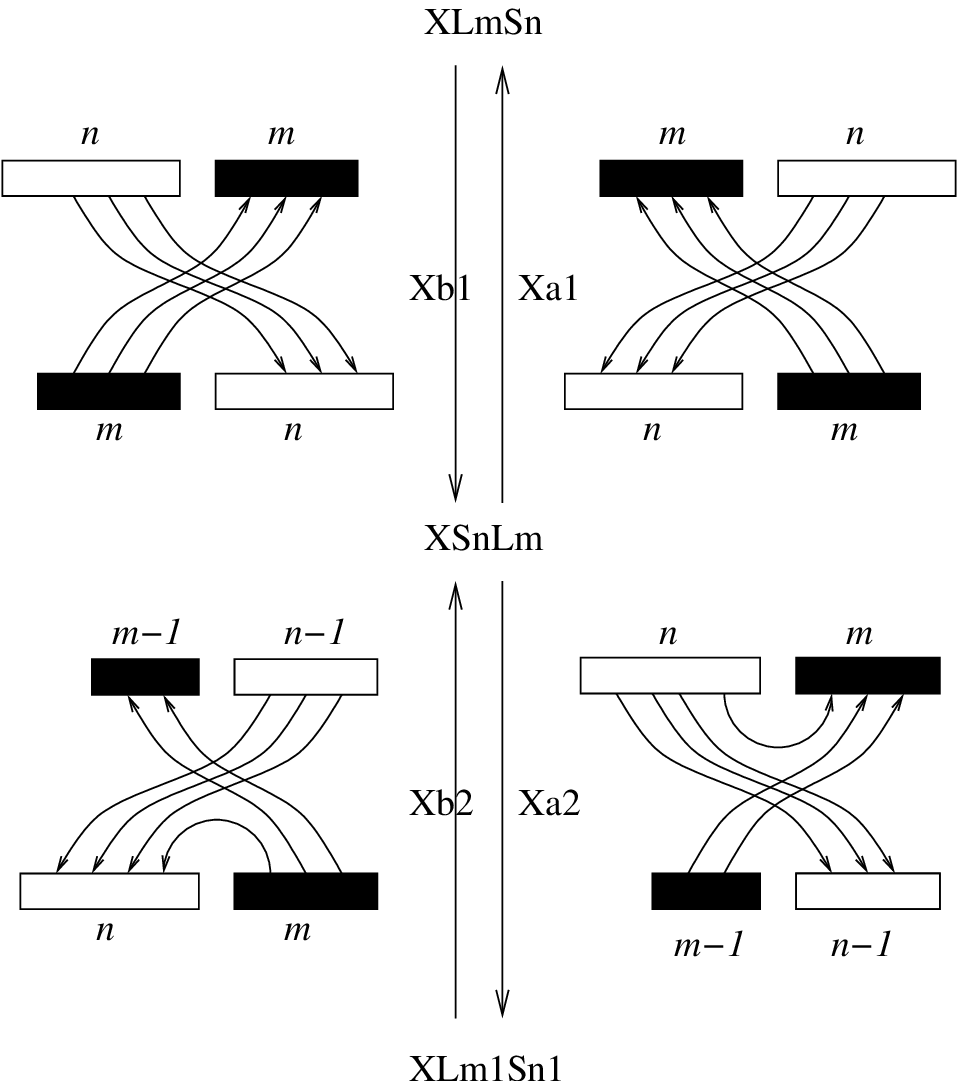}}\end{center}

A straightforward manipulation of diagrams shows that 
$$  \alpha_1 \beta_1 = \Id, \ \  \alpha_1\beta_2 =0, \ \ 
 \alpha_2\beta_2 = \frac{1}{mn}  \Id, \ \ \alpha_2 \beta_1 =0.  $$
Let $\beta'_2 = mn \beta_2$. Then the maps 
$$\Lambda^m_+\otimes S^n_-  
\displaystyle_{\stackrel{\longrightarrow}{\beta_1}}^{\stackrel{\alpha_1}{\longleftarrow}}
 S^n_- \otimes \Lambda^m_+ 
\displaystyle_{\stackrel{\longleftarrow}{\beta_2'}}^{\stackrel{\alpha_2}{\longrightarrow}}
  \Lambda^{m-1}_+ \otimes S^{n-1}_- $$
satisfy
$$ \alpha_1 \beta_1 = \Id, \ \  \alpha_1\beta'_2 =0,
\ \  \alpha_2\beta'_2 = \Id, \ \ \alpha_2 \beta_1 =0, \ \  
   \beta_1\alpha_1 + \beta'_2 \alpha_2 = \Id. $$ 
The last equality follows from a direct diagrammatic manipulation as well. 
Thus, there is an isomorphism 
$$ S^n_- \otimes \Lambda^m_+  \cong \big(\Lambda^m_+\otimes S^n_- 
\big) \oplus \big( \Lambda^{m-1}_+ \otimes S^{n-1}_- \big) ,$$ 
concluding the proof of Proposition~\ref{prop-can-intr} and 
Corollary~\ref{cor-rels}.

There is a natural bijection between partitions $\lambda$ of $n$ and 
(isomorphism classes of) irreducible representations of $\Bbbk[S_n]$. 
To each partition $\lambda=(\lambda_1, \dots, \lambda_k)$ with 
$|\lambda|=\lambda_1+\dots + \lambda_k=n$ 
there corresponds the unique common irreducible summand $L_{\lambda}$ 
of the representation induced from the trivial representation of parabolic subgroup 
$S_{\lambda} = S_{\lambda_1}\times \dots \times S_{\lambda_k}\subset S_n$ 
and the representation induced from the sign representation of parabolic subgroup 
$S_{\lambda^{\ast}}= S_{\lambda^{\ast}_1}\times \dots \times
 S_{\lambda^{\ast}_m}$, where $\lambda^{\ast}$ is the dual partition. 
Let $e_{\lambda}\in \Bbbk[S_n]$ be the Young idempotent, so that 
$e_{\lambda}^2= e_{\lambda}$ and $L_{\lambda}\cong \Bbbk[S_n]e_{\lambda}.$ 

We denote by $Q_{+,\lambda}:= (Q_{+^n}, e_{\lambda})$ the 
object of $\Ctwo$ which is the direct summand of $Q_{+^n}$ corresponding to the
idempotent $e_{\lambda}$, where we view the latter as an idempotent in the endomorphism 
ring via the standard homomorphism  $\Bbbk[S_n]\lra \End_{\Ctwo}(Q_{+^n})$. 
Likewise, let 
$Q_{-,\lambda}:= (Q_{-^n}, e_{\lambda})$ be the corresponding direct summand 
of $Q_{-^n}$, where we view $e_{\lambda}$ as an endomorphism of the latter 
object. In particular, 
$$ S^n_+ = Q_{+,(n)}, \ \  \Lambda^n_+ = Q_{+, (1^n)}, \ \ 
   S^n_- = Q_{-,(n)}, \ \  \Lambda^n_- = Q_{-,(1^n)}.$$ 

\vspace{0.06in} 

The Grothendieck ring $K_0(\Ctwo)$ is an abelian group with generators--symbols 
$[M]$, over all objects $M$ of $\Ctwo$ and defining relations 
$[M_1]=[M_2]+[M_3]$ whenever $M_1\cong M_2 \oplus M_3$. 
Monoidal structure on $\Ctwo$ descends to an associative multiplication 
on $K_0(\Ctwo)$, with $[\mathbf{1}]$ being the identity for multiplication. 
Hence, $K_0(\Ctwo)$ is an associative unital ring. 

\vspace{0.06in} 

Recall the ring $H_{\Z}$ from the introduction. We can now define homomorphism 
$\gamma : H_{\Z}\lra K_0(\Ctwo)$ discussed there:  
$$ \gamma(a_n) = [Q_{-,(n)}]=[S^n_-], \ \ \ \ 
\gamma(b_m)=[Q_{+, (1^m)}]= [\Lambda_+^m].$$ 

If we identify the subring of $H_{\Z}$ generated by the $a_n$'s 
with the ring of symmetric functions $\mathrm{Sym}$ so that $a_n$ 
corresponds to $n$-th complete symmetric function $h_n$, then 
$\gamma$ will take the Schur function associated to partition $\lambda$ 
to $[Q_{+,\lambda}]$. This function is often denoted $s_{\lambda}$; for us 
it is convenient to call it $a_{\lambda}$, so that $a_{(n)}=a_n$. 

Similarly,  we identify the subring generated by the $b_m$'s with 
$\mathrm{Sym}$ by taking $b_m$ to the $m$-th elementary symmetric function 
$e_m$.  Denote by $b_{\lambda}$ the polynomial in $b_m$'s that corresponds to 
the Schur function $s_{\lambda}$ under this identification. In particular,  
$b_{(1^m)}=b_m$. We have 
$$ \gamma(a_{\lambda}) =  [Q_{-, \lambda}], \ \ \ \ 
\gamma(b_{\lambda}) =  [Q_{+, \lambda^{\ast}}]. $$

Littlewood-Richardson coefficients $r_{\lambda,\mu}^{\nu}$ that 
appear in decompositions of the product of Schur functions 
\begin{eqnarray*} 
a_{\lambda} a_{\mu} & = & \sum_{\nu} r_{\lambda,\mu}^{\nu} a_{\nu}, \\
b_{\lambda} b_{\mu} & = & \sum_{\nu} r_{\lambda,\mu}^{\nu} b_{\nu},
\end{eqnarray*} 
also appear in the isomorphism formulas in $\Ctwo$: 
\begin{eqnarray*} 
Q_{+,\lambda}\otimes Q_{+,\mu}  & \cong & 
\oplusop{\nu} (Q_{+,\nu})^{r_{\lambda,\mu}^{\nu}} , \\
Q_{-,\lambda}\otimes Q_{-,\mu}  & \cong & 
\oplusop{\nu} (Q_{-,\nu})^{r_{\lambda,\mu}^{\nu}} . 
\end{eqnarray*} 
Descending to the Grothendieck ring, we have 
\begin{eqnarray*} 
 \  [Q_{+,\lambda}][Q_{+,\mu}]  & = & 
\sum_{\nu} r_{\lambda,\mu}^{\nu} [Q_{+,\nu}] ,  \\     
   \    [    Q_{-,\lambda} ] [Q_{-,\mu}]  & = & 
\sum_{\nu} r_{\lambda,\mu}^{\nu} [Q_{-,\nu}] . 
\end{eqnarray*} 

The ring $H_{\Z}$ has a basis $\{b_{\mu}a_{\lambda}\}_{\lambda,\mu}$ over 
all partitions $\lambda, \mu$. Consequently, elements 
$ [Q_{+, \mu}][Q_{-,\lambda}]$ over all $\lambda,\mu$ span the subring 
$\gamma(H_{\Z})$ of $K_0(\Ctwo)$. 

\vspace{0.06in} 

\emph{Remark:} 
The symmetries $\xi_1, \xi_2, \xi_3$ of 
$\Cone$ extend to self-equivalences of category $\Ctwo$, also denoted 
$\xi_1,$ $\xi_2,$ $\xi_3$. On objects $Q_{+,\lambda}\otimes Q_{-,\mu}$ 
they act as follows: 
\begin{eqnarray*} 
 \xi_1 (Q_{+,\mu}\otimes Q_{-,\lambda}) & = & Q_{+,\mu^{\ast}}\otimes 
Q_{-,\lambda^{\ast}}, \\
\xi_2 (Q_{+,\mu}\otimes Q_{-,\lambda}) & = & Q_{+,\mu}\otimes Q_{-,\lambda}, \\
\xi_3(Q_{+,\mu}\otimes Q_{-,\lambda}) & = &Q_{+,\lambda}\otimes Q_{-,\mu}. 
\end{eqnarray*} 
These self-equivalences induce involutions $[\xi_1], [\xi_2]$ and 
antiinvolution $[\xi_3]$ on $K_0(\Ctwo)$.  The involution of $H_{\Z}$ corresponding 
to $[\xi_2]$ is the identity. We don't know whether $[\xi_2]$ is the identity 
involution on the entire $K_0(\Ctwo)$; this would follow from Conjecture~\ref{conj-iso}.

%%%%%%%%%%%%%%%%%%%%%%%%%%%
%%%%%%%%%%%%%%%%%%%%%%%%%%%
%%
%% INDUCTION AND RESTRICTION, FINITE GROUPS
%%
%%%%%%%%%%%%%%%%%%%%%%%%%%%
%%%%%%%%%%%%%%%%%%%%%%%%%%%

\section{Diagrammatics for induction and restriction functors} 
\label{sec-indres}

%%%%%%%%%%%%%%%%%%%%%%%%%%
%%
%%   BIADJOINT FUNCTORS 
%%
%%%%%%%%%%%%%%%%%%%%%%%%%%

\subsection{Biadjoint functors} \label{subsec-biadjoint}

Recall~\cite{MacLane} that a functor $L: \mc{A}\lra \mc{B}$ between categories 
$\mc{A}$ and $\mc{B}$ is left adjoint to a functor $R: \mc{B} \lra \mc{A}$ 
whenever there are natural transformations 
\begin{equation}\label{eq-adj-maps1}
\alpha: LR \Rightarrow \Id_{\mc{B}}, \ \ \ 
 \beta: \Id_{\mc{A}}\Rightarrow RL 
\end{equation}
that satisfy the relations 
\begin{equation}\label{eq-biadj1}
 (\alpha \circ \Id_L )(\Id_L \circ \beta) = \Id_L, \ \ 
  (\Id_R \circ \alpha)  (\beta \circ \Id_R)= \Id_R.
\end{equation} 

Assume that $L$ is both left and right adjoint to $R$, and the second 
adjunction maps 
\begin{equation}\label{eq-adj-maps2}
\overline{\alpha}: RL \Rightarrow \Id_{\mc{A}}, \ \ \ 
 \overline{\beta}: \Id_{\mc{B}}\Rightarrow LR 
\end{equation} 
are fixed as well. They satisfy     
\begin{equation}\label{eq-biadj2}
(\overline{\alpha} \circ \Id_R) (\Id_R \circ {}\overline{\beta})= \Id_R, \ \ 
(\Id_L \circ \overline{\alpha}) (\overline{\beta} \circ \Id_L) = \Id_L.
\end{equation}

Out of $\alpha,\overline{\alpha},\beta,\overline{\beta}$ one can construct 
more general natural transformations between compositions of functors 
$L$ and $R$ by placing the basic four transformation in various locations 
in the composition of functors, and then composing several such transformations. 
It's convenient to draw these compositions via planar diagrams, 
with transformations (\ref{eq-adj-maps1}), (\ref{eq-adj-maps2}) depicted as U-turns 

\begin{center}{\psfig{figure=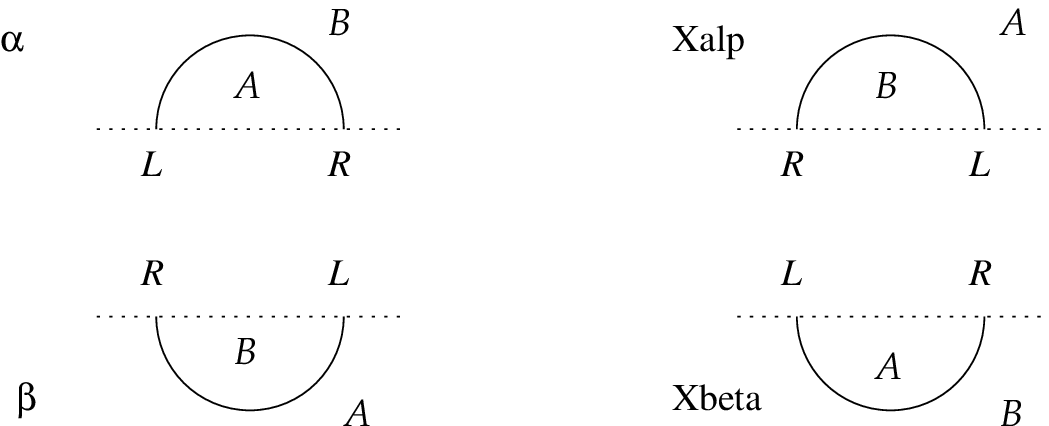}}\end{center}

General diagrams are built out of U-turns and vertical lines, the latter 
denoting identity natural transformations of $R$ and $L$. For instance, 

\begin{center}{\psfig{figure=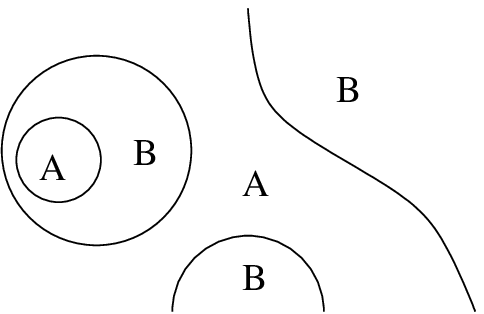}}\end{center}

is the following natural transformation from $RLR $ to $R$ 

$$(\overline{\alpha}\circ \Id)  (\Id \circ \alpha\circ \Id^{\otimes 2}) 
(\Id \circ \overline{\beta}\circ \Id^{\otimes 2})   
(\beta\circ \Id)  \overline{\alpha} $$

The four biadjointness equations (\ref{eq-biadj1}), (\ref{eq-biadj2}), which can 
drawn as   

\begin{center}{\psfig{figure=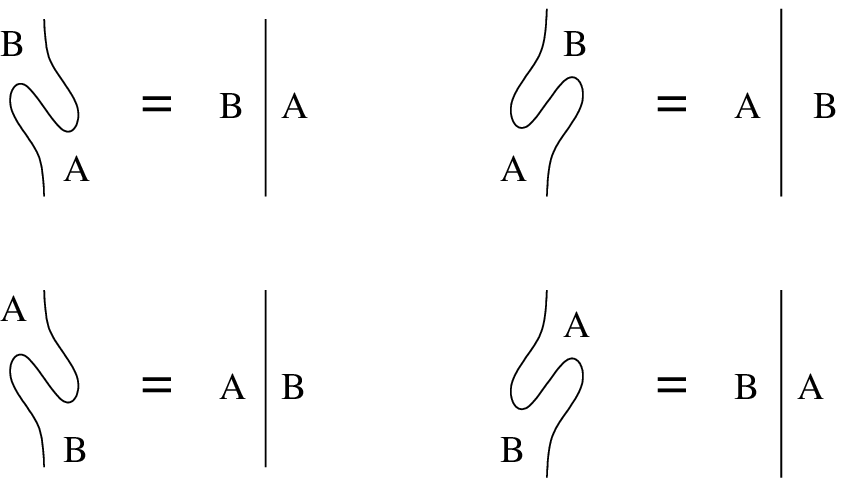}}\end{center}

\noindent 
are equivalent to the condition that the lines and circles can be isotoped without 
changing the natural transformation associated to the diagram. This was 
observed in~\cite{CKS, Muger}. 

The graphical calculus of biadjoints can be further enhanced. Assume given 
a collection of categories and a collection of functors between them such 
that each functor has a biadjoint, which is also in the collection, and the 
biadjointness transformations are fixed. Natural transformations 
generated by the biadjointness ones can be drawn via diagrams on the 
plane strip $\R\times [0,1]$, with lines and circles labelled by functors 
and regions labelled by categories, with arbitrary rel boundary isotopies 
allowed. 

Furthermore, any element $z$ in the center of a category $\mc{A}$
(i.e. $z$ is the endomorphism of the identity functor $\Id_{\mc{A}}$) 
can be shown as freely floating in a region labelled $\mc{A}$. 
Two such central elements can freely move past each other. 

Given two functors $F,G: \mc{A} \lra \mc{B}$  in the collection, a natural 
transformation $a: F\lra G$ can be depicted as a labelled dot 
on a line separating a segment labelled $F$ from a segment labelled $G$. 

\begin{center}{\psfig{figure=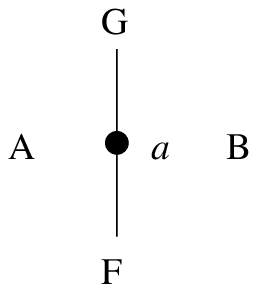}}\end{center}

When $a$ is dragged through a $U$-turn, it can change in two 
possible ways, depending on the type of a $U$-turn, into natural transformations 
$a^{\ast}, {}^{\ast}a: G'\lra F'$, where $F',G'$ are the biadjoints of 
$F$ and $G$. If ${}^{\ast} a=a^{\ast}$, we say that the 
biadjointness data is \emph{cyclic}. For more on the cyclicity condition 
see~\cite{Bartlett1, Bartlett2, CKS, AL1}.  

Biadjoint functors appear throughout categorification and TQFTs, 
see discussions in~\cite[Section 5]{FVIT}, \cite{MK} and references therein, 
also~\cite{Brun}.

%%%%%%%%%%%%%%%%%%%%%%%%%
%%
%% IND AND RES BETWEEN FINITE GROUPS
%%
%%%%%%%%%%%%%%%%%%%%%%%%%

\subsection{Induction and restriction between finite groups in pictures}
\label{subsec-indres-pic}

We fix a commutative ring $\Bbbk$ and denote by  
$\Bbbk G$ the group algebra of  finite group $G$ with coefficients in $\Bbbk$. 
We denote by $(G)$ the group algebra $\Bbbk G$, viewed as a 
$(\Bbbk G, \Bbbk G)$-bimodule. 
If $H$ is a subgroup of $G$ and we view $\Bbbk G$ as a $(\Bbbk G,\Bbbk H )$-bimodule 
via the left action of $\Bbbk G$ and the right action of $\Bbbk H$, 
we denote it by $(G_H).$ When viewing  $\Bbbk G$ as a $\Bbbk H$-bimodule, 
denote it by $({}_HG_H)$, etc. Similar shortcut notation is adopted for 
tensor products of bimodules. For instance, $(G_H G)$ denotes 
$\Bbbk G$-bimodule $\Bbbk G\otimes_{\Bbbk H}\Bbbk G$, while 
$({}_H G_H G)$  denotes the same space, but viewed as a $(\Bbbk H,\Bbbk G)$-bimodule.  

Start with the 2-category $\finone$ whose objects are finite groups $G$, 
morphisms from $G$ to $H$ are $(\Bbbk H,\Bbbk G)$-bimodules, and 
2-morphisms are bimodule homomorphisms. 
Consider the 2-subcategory $\fintwo$ of $\finone$ 
with the same objects as $\finone,$ while morphisms are finite direct 
sums of tensor products of bimodules $(G_H)$ and $(_H G)$ corresponding 
to the induction and restriction functors between categories of $H$ and $G$-modules. 
Thus, a 1-morphism from $G$ to $G'$ is a finite direct sum of bimodules 
isomorphic to 
$$ (G_n \ {}_{H_{n-1}} G_{n-1} \ {}_{H_{n-2}} \dots {}_{H_2} G_2 \ {}_{H_1} G_1) $$ 
where $G'=G_n \supset H_{n-1}\subset G_{n-1}\supset \dots \subset G_2 \supset H_1\subset 
G_1=G$ is a zigzag of inclusions between finite groups. 
The 2-morphisms in $\fintwo$ are bimodule homomorphisms. 
Alternatively, we can think of 1-morphisms in $\fintwo$ as given by compositions 
of induction and restriction functors between categories of (left) $G$-modules, 
over finite groups $G$. 

In this section, we develop basics of a  graphical calculus for studying 
2-morphisms in $\fintwo$. 
In general, given a unital inclusion of rings $B\subset A$, the induction 
functor $\Ind: B \dmod\lra A\dmod$ that takes $M$ to $A\otimes_B M$ is 
left adjoint to the restriction functor. An inclusion $\iota: H\subset G$ of 
finite groups produces an inclusion $\Bbbk H\subset \Bbbk G$ of group algebras, 
with the induction functor 
$$ \Ind_H^G \ : \  \Bbbk H\dmod \lra \Bbbk G\dmod$$ 
being both left and right adjoint (i.e. biadjoint) 
to the restriction functor 
$$ \Res_G^H \ : \  \Bbbk G\dmod \lra \Bbbk H\dmod.$$ 
The biadjointness endomorphisms are given by the following four bimodule maps 

\vspace{0.06in} 

\noindent 
1) $ \  : \ (G_H G ) \lra (G), \ \ x\otimes y \longmapsto xy, \ \  x,y \in (G), $ 

\vspace{0.06in} 

\noindent 
2)  $\ : \  (H) \lra ({}_H G _G G_H ) , \ \  x \longmapsto x\otimes 1 = 1 \otimes x, 
 \ \  x\in (H) , $

\vspace{0.06in} 

\noindent
$ 3)  \ : \ ({}_H G_G G_H ) \cong ({}_HG_H) \lra (H), \ \ 
   g \longmapsto g \mbox{ if }g\in H, 
   g \longmapsto 0, \mbox{ if }g \in G\setminus H. $ 

\vspace{0.06in} 

\noindent 
We denote this projection map by  $p_H: ({}_HG_H) \lra (H)$, 
$p_H(g)=g$ if $g\in H$ and $p_H(g)=0$ if $g\in G\setminus H$, 
extended by $\Bbbk$-linearity. Clearly, $p_H$ is a map of $\Bbbk H$-bimodules. 

\vspace{0.06in} 

\noindent 
4) Let $G=\bigsqcup_{i=1}^m Hg_i$ be a decomposition of $G$ into left $H$-cosets, 
so that $m=[G:H]$ is the index of $H$ in $G$. Notice that the element 
$$\sum_{i=1}^m g_i^{-1} \otimes g_i  \ \in (G_H G)$$
does not depend on the choice of coset representatives $\{g_i\}_{i=1}^m$ 
of $H$ in $G$: if $g'_i=h_i g_i$ then 
$$ \sum_{i=1}^m g'^{-1}_i \otimes g'_i =   \sum_{i=1}^m g^{-1}_ih_i^{-1} 
\otimes h_i g_i = \sum_{i=1}^m g_i^{-1} \otimes g_i,$$ 
since the tensor product is over $\Bbbk[H]$, and $h_i$ can be moved 
through the tensor product sign. Define bimodule map 
\begin{equation}\label{eq-dcup} (G) \lra (G_H G)  \end{equation} 
by the condition that 
$$ 1 \longmapsto \sum_{i=1}^m g_i^{-1} \otimes g_i ,$$
so that 
$$ g \longmapsto \sum_{i=1}^m g_i^{-1} \otimes g_i g = 
\sum_{i=1}^m g g_i^{-1} \otimes g_i.$$ 
The second equality, needed to insure that one does get a bimodule map, 
follows from the following computation: $g_i g= h_i g_{i'}$ for some
 $h_i\in H$ and $i'\in \{1,2,\dots, m\}$. The assigment $i\longmapsto i'$ 
is a bijection of $\{1,2,\dots, m\}$. We have 
\begin{equation*} 
\sum_{i=1}^m g_i^{-1} \otimes g_i  g  =
\sum_{i=1}^m g_i^{-1} \otimes h_i g_{i'}   =  
 \sum_{i=1}^m g_i^{-1}h_i  \otimes g_{i'}   = 
\sum_{i'=1}^m g g_{i'}^{-1}  \otimes g_{i'} = 
\sum_{i=1}^m g g_i^{-1} \otimes g_i. 
 \end{equation*}

\noindent 
Combining with the earlier remark, we see that (\ref{eq-dcup})  is a bimodule map 
which does not depend on the choices of $H$-coset representatives $g_i$. 

\vspace{0.06in} 

We associate to these four bimodule maps the following four pictures

\vspace{0.15in} 

{\psfig{figure=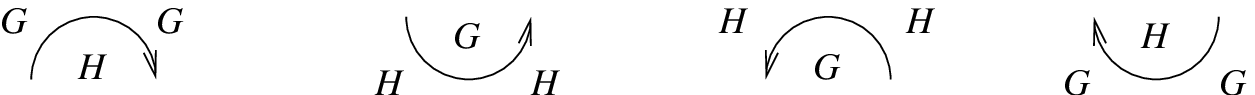}}

\vspace{0.15in} 

Thus, 

\raisebox{-0.4cm}
{\psfig{figure=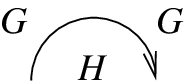}}, denoted $\alpha_H^G$, is the map 
$$ (G_H G ) \lra (G), \ \ x\otimes y \longmapsto xy, \ \  x,y \in (G).$$

\raisebox{-0.4cm}
{\psfig{figure=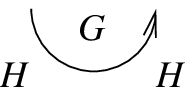}}, denoted $\beta_H^G$, is the map 
$$(H) \lra ({}_H G _G G_H ) \cong ({}_H G_H), \ \ x\longmapsto x\otimes 1 = 
1\otimes x , \ \  x\in (H).$$ 

\raisebox{-0.4cm}
{\psfig{figure=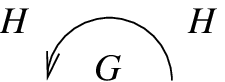}}, denoted $\overline{\alpha}_G^H$, is the map $p_H$ 
described earlier,  
$$ ({}_H G_G G_H ) \cong ({}_HG_H) \lra (H), \ \  x \longmapsto p_H(x), 
\ \  x\in (G).$$ 

\raisebox{-0.4cm}
{\psfig{figure=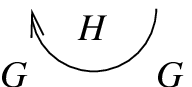}}, denoted $\overline{\beta}_G^H$, is the $\Bbbk$-linear map 
(\ref{eq-dcup}) 
$$  (G) \lra (G_H G), \  \  g \longmapsto \sum_{i=1}^m g_i^{-1}\otimes g_i g , \ \
 g\in G.$$

\begin{theorem} \label{theorem-cbp}
These four bimodule maps turn induction and restriction 
functors $\Ind_H^G$ and $\Res_G^H$ into a cyclic biadjoint pair. 
\end{theorem} 

\emph{Proof:} 
First, we check that the adjointness equations (\ref{eq-biadj1}) and (\ref{eq-biadj2}) 
hold for these maps. 
The bimodule map $ (G_H) \lra (G_H G_H) \lra (G_H)$ corresponding
to the left hand side of the first equation in (\ref{eq-biadj1}) is given by 
$g \longmapsto g \otimes 1 \longmapsto g1=g$, hence the map is 
the identity: 

\begin{center}{\psfig{figure=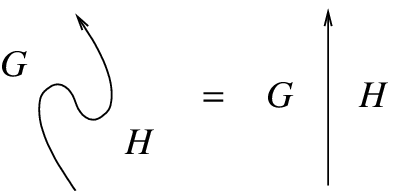}}\end{center}

The bimodule map $({}_HG) \lra ({}_HG_H G) \lra ({}_HG)$ for the 
left hand side of the second equation in  (\ref{eq-biadj1}) is given by 
$g \longmapsto 1\otimes g \longmapsto 1g= g$, and the map is 
the identity: 

\begin{center}{\psfig{figure=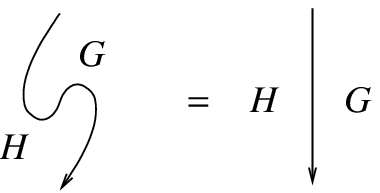}} \end{center} 
 
The bimodule map $({}_HG) \lra ({}_HG_H G) \lra ({}_HG)$ for the 
left hand side of the first equation in  (\ref{eq-biadj2}) is given by 
$$g \longmapsto
\sum _{i=1}^m gg_i^{-1}\otimes g_i\longmapsto 
\sum _{i=1}^m p_H(gg_i^{-1}) g_i .$$ 
Notice that $p_H(gg_i^{-1})=0$ iff $g\notin Hg_i$, and 
$p_H(hg_i g_i^{-1})=h$. Therefore, $g\longmapsto g$ under the map,  
and the first equation in  (\ref{eq-biadj2}) holds:  

\begin{center}{\psfig{figure=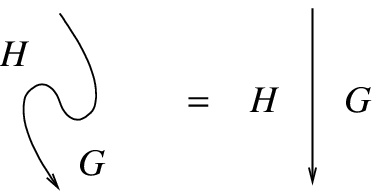}}\end{center} 
 
The bimodule map $(G_H) \lra (G_HG_H) \lra (G_H)$ for the 
left hand side of the second equation in  (\ref{eq-biadj2}) is given by 
$$g\longmapsto \sum _{i=1}^m g_i^{-1}\otimes g_ig \longmapsto 
\sum_{i=1}^m g_i^{-1} p(g_ig)=g$$ 
by a similar computation, so that  

\begin{center}{\psfig{figure=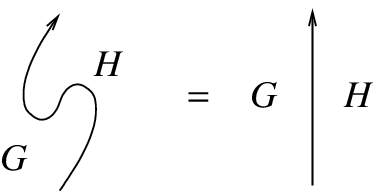}}, \end{center} 

\noindent 
and the four bimodule maps above determine biadjointness morphisms for induction 
and restriction functors between $H$ and $G$. 

Consider the $\Bbbk$-algebra 
$$ \Bbbk G^H \ := \ \{ a\in \Bbbk G| ha=ah \ \  \forall h\in H \} $$ 
of $H$-invariants in $\Bbbk G$ with respect to the conjugation action. 
This algebra is canonically isomorphic to the endomorphism ring 
of the bimodule $({}_HG)$, and, therefore, to the endomorphism ring 
of the functor $\Res_G^H$, via the map that assigns to $a\in \Bbbk G^H$ 
the endomorphism ${}'{a}(x) := a x $, where $x\in ({}_HG)$. 
Likewise, the opposite algebra of $\Bbbk G^H$ is canonically isomorphic 
to the endomorphism ring of the bimodule $(G_H)$ and, therefore, to that of  
the functor $\Ind_H^G$, via the map that assigns to $a\in \Bbbk G^H$ 
the endomorphism $a'(x):= xa$, where $x\in (G_H)$. 

Thus, to $a\in \Bbbk G^H$ we assign endomorphisms $a'$ and ${}'a$ 
of $\Ind_H^G$ and $\Res_G^H$ and depict them by 

\begin{center}\psfig{figure=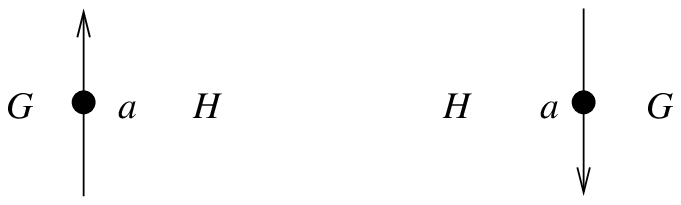}\end{center}

\begin{lemma} \label{lemma-4eq} 
For any $a\in  \Bbbk G^H$ there are equalities 
of bimodule homomorphisms 
\begin{eqnarray} 
 \raisebox{-0.6cm}{\psfig{figure=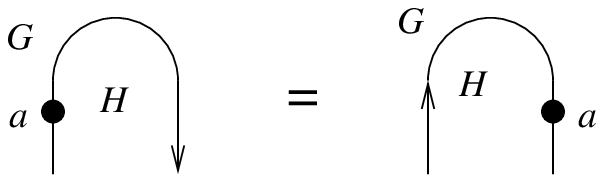}}  & &  \alpha_H^G (a' \circ \Id)  = 
\alpha_H^G (\Id \circ {}'a), \\
\raisebox{-0.6cm}{\psfig{figure=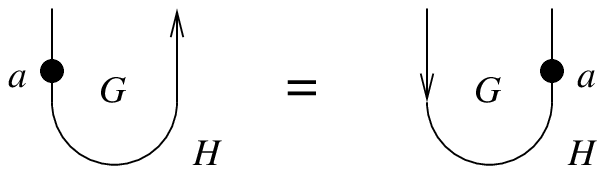}} & & 
 ({}'a\circ \Id) \beta_H^G   =    (\Id \circ a') \beta_H^G, \\
\raisebox{-0.6cm}{\psfig{figure=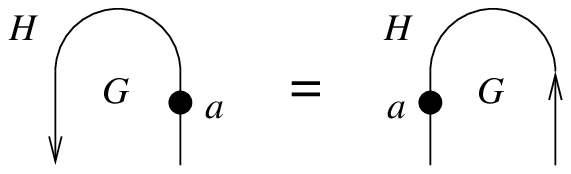}} & & 
 \overline{\alpha}_G^H ({}'a\circ \Id)  =  \overline{\alpha}_G^H 
  (\Id \circ a'), \\
\raisebox{-0.6cm}{\psfig{figure=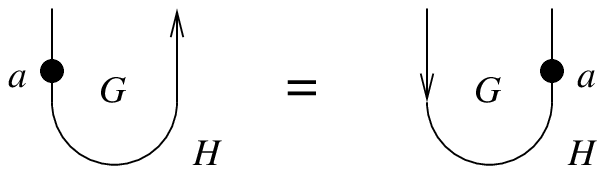}} & & 
   (a' \circ \Id ) \overline{\beta}_G^H  =  (\Id \circ {}'a) \overline{\beta}_G^H . 
\end{eqnarray}
\end{lemma}

The left hand side of the first equality is a map $(G_H G)\lra (G)$ given by 
 $g_1\otimes g_2 \longmapsto g_1 a \otimes g_2 \longmapsto g_1 ag_2 .$ 
The right hand side is 
$ g_1 \otimes g_2 \longmapsto g_1 \otimes ag_2 \longmapsto g_1 ag_2,$ 
and the equality is obvious. 

The second equality follows from an equally trivial computation. 

The third equality is the equation 
$ p_H(ga) = p_H(ag) $ for $g\in G$ and $a\in \Bbbk G^H$. It suffices 
to check it when $\Bbbk = \Z$ and $a=\sum_{h\in H} h k h^{-1}$ for 
some $k\in G$. The equation becomes 
\begin{equation} \label{eq-ghkh}
 \sum_{h \in H} p_H(g h k h^{-1}) = \sum_{h\in H} p_H(hkh^{-1}g) 
\end{equation} 
The left hand side equals
$$
\sum_{h \in H} p_H(g h k) h^{-1} = \sum_{h^{-1}g^{-1}u=k} p_H(u)p_H(h^{-1}) 
 =   \sum_{hg^{-1}u=k} p_H(u)p_H(h)$$
where in the first equality we set $u= ghk$, the sum being over all $u,h\in G$ with 
$h^{-1}g^{-1}u=k$. For the second equality, we converted $h$ to $h^{-1}$. 

The right hand side of (\ref{eq-ghkh}) equals 
$$  \sum_{h\in H} h p_H(kh^{-1}g) = \sum_{u g^{-1}h=k}p_H(h)p_H(u), $$ 
where we set $u=k h^{-1}g$ and the sum is over all $u,h\in G$. Interchanging 
$h$ and $u$, we see that (\ref{eq-ghkh}) holds. 

For the last of the four equations, 
it suffices to check that the image of $1\in G$ is the same under these two 
bimodule homomorphisms. For the one on the left, 
$$ 1\longmapsto \sum_{i=1}^m g_i^{-1} \otimes g_i \longmapsto 
 \sum_{i=1}^m g_i^{-1}\otimes a g_i . $$
For the one on the right, 
$$ 1 \longmapsto \sum_{i=1}^m g_i^{-1} \otimes g_i \longmapsto 
\sum_{i=1}^m g_i^{-1} a \otimes g_i . $$
Again, we can assume $\Bbbk=\Z$ and $a=\sum_{h\in H} h k h^{-1}$ 
for some $k\in G$. The equation becomes 
\begin{equation} \label{eq-sumih}
\sum_{i, h} g_i^{-1}\otimes h k h^{-1} g_i  = 
 \sum_{i, h} g_i^{-1} h k h^{-1} \otimes g_i, 
\end{equation} 
summing over $1\le i \le m$ and $h\in H$. We have 
$$ \sum_{i, h} g_i^{-1}\otimes h k h^{-1} g_i = 
 \sum_{i, h} g_i^{-1}h \otimes k h^{-1} g_i = 
\sum_{u\in G} u \otimes k u^{-1}, $$ 
where, in the first equality, $h$ is moved to the left (the tensor 
product is over $\Bbbk H$), and in the second equality $u= g_i^{-1}h$ 
runs over all elements of $G$ as $i$ changes from $1$ to $m$ and 
$h$ runs over all elements of $H$. Likewise, 
$$ \sum_{i, h} g_i^{-1} h k h^{-1} \otimes g_i = 
  \sum_{i, h} g_i^{-1} h k \otimes h^{-1} g_i = 
  \sum_{u\in G} u k \otimes u^{-1}.$$ 
Equation (\ref{eq-sumih}) and Lemma~\ref{lemma-4eq} follow.

Since the box labelled by $a\in \Bbbk G^H$ can be dragged through 
any $U$-turn, we see that the biadjointness maps have the cyclic property
-- dragging the box labelled $a$ all the way along a circle brings us 
back to the original diagram. 
This concludes the proof of Theorem~\ref{theorem-cbp}. $\square$

\vspace{0.1in} 

There are obvious simplification relations 

\begin{center}{\psfig{figure=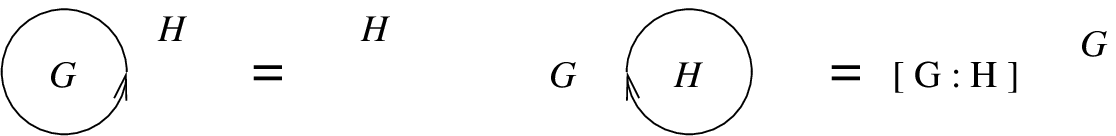}}\end{center}
 
The first relation says that a counterclockwise bubble with $G$ inside 
and $H$ outside can be erased. The second relation allows to remove 
a clockwise bubble at the cost of multiplying the diagram by the index 
of $H$ in $G$. 

\vspace{0.1in} 

For each inclusion of finite groups $H\subset K \subset G$ 
there is a canonical isomorphism between induction functors 
$\Ind^G_H\cong \Ind^G_K\circ \Ind^K_H$ which corresponds 
to a canonical isomorphism of bimodules 
 $(G)_H\cong (G)_K (K)_H$. Likewise, canonical isomorphism 
between restrictions $\Res_G^H \cong \Res_K^H \circ \Res_G^K$
is given by a natural isomorphism of bimodules 
${}_H(G) \cong {}_H(K)_K(G)$. 
We draw these isomorphisms via trivalent diagrams  

\begin{center}{\psfig{figure=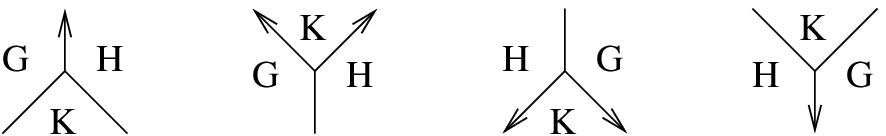}}\end{center}
 
Since the isomorphisms are mutually-inverse, we have, for the 
first two isomorphisms, 
\begin{center}{\psfig{figure=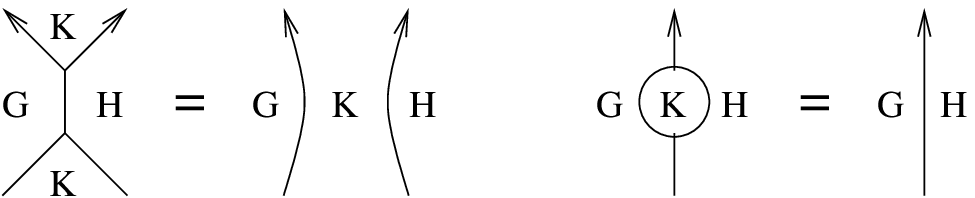}}\end{center}

Mutual inversion of the other two isomorphisms can be similarly depicted.  
These definitions are compatible with isotopies--identities

\begin{center}{\psfig{figure=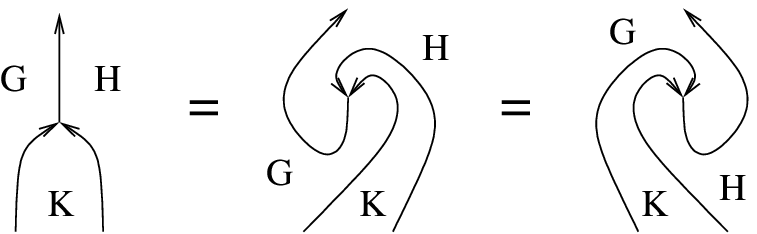}}\end{center}

\noindent 
hold (likewise for the other pair of isomorphisms). These identities 
imply that various definitions of trivalent vertices in other positions 
relative to the $y$-coordinate are all the same. For instance, 
if we define 

\begin{center}{\psfig{figure=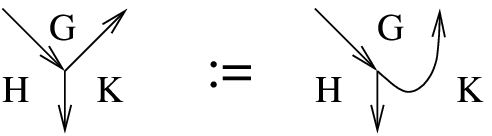}}\end{center}

then 

\begin{center}{\psfig{figure=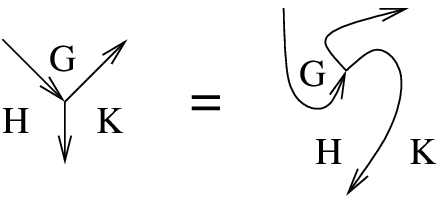}}\end{center}

\noindent
The associativity relation for the  
inclusions of four groups $L\subset H \subset K \subset G$ has the form 

\begin{center}{\psfig{figure=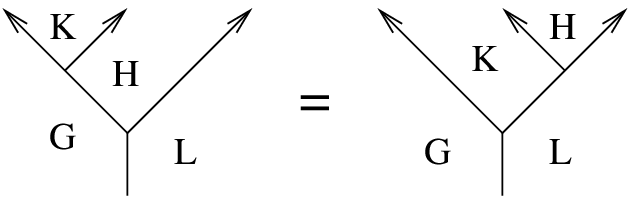}}\end{center}
 
The induction and restriction functors for $H$ and $G$ depend on 
the inclusion $\iota: H\hookrightarrow G$. Conjugating the inclusion by 
an element $g\in G,$ so that $\iota'(h)= g h g^{-1},$ $\iota' : H\hookrightarrow G$,
 leads to induction and restriction functors isomorphic to the 
original ones, via bimodule maps 
\begin{eqnarray*} 
& & (G)_{\iota(H)} \lra  (G)_{\iota'(H)}, \ \  f \longmapsto fg^{-1} , \ \  f\in G, \\
& & {}_{\iota(H)}(G) \lra {}_{\iota'(H)}(G), \ \  f \longmapsto gf, \ \ f\in G.
\end{eqnarray*} 
We depict these conjugation isomorphisms via a mark on a line with 
$g$ next to it: 
\begin{center}{\psfig{figure=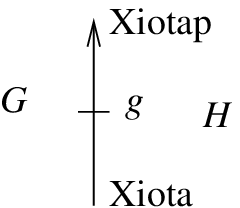}}\end{center}

The Mackey induction-restriction theorem says that, given 
subgroups $H,K$ of a finite group $G$, there is an isomorphism 
\begin{equation} 
\Res_G^K \circ \Ind^G_H \cong \oplusop{i\in I} \Ind^K_{K\cap g_i H g_i^{-1}} 
\circ \Res_H^{K\cap g_i H g_i^{-1}} ,
\end{equation} 
where the sum is over representatives $g_i$ of  $(K,H)$-cosets of $G$, 
$$  G = \sqcup_{i\in I} K g_i H.$$ Let $L_i = K\cap g_i H g_i^{-1}$. 
Diagrams 
\begin{center}{\psfig{figure=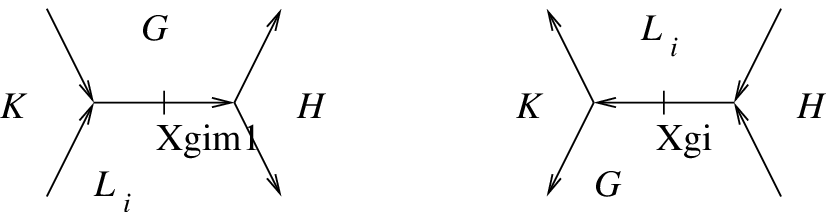}}\end{center}
\noindent 
define $(K,H)$-bimodule maps 
$$ \alpha_i \ : \ (K)_{L_i}{}'(H) \lra {}_K(G)_H, \ \ \ \beta_i \ : \ {}_K(G)_H \lra 
  (K)_{L_i}{}'(H).$$ 
Here $(K)_{L_i}{}'(H)$ is $\Bbbk[K]\otimes_{\Bbbk[L_i]}\Bbbk[H]$, with 
$x\in L_i$ acting on $H$ by right multiplication by $g_i^{-1}xg_i$. 

\begin{prop} The maps $\sum_{i\in I} \alpha_i$ and $\sum_{i\in I} \beta_i$ are 
mutually-inverse isomorphisms of bimodules 
 $\oplus_{i\in I} (K)_{L_i}{}'(H)$  and ${}_K(G)_H.$ 
\end{prop} 

This proposition is a pictorial restatement of the Mackey theorem. 
Proof is left to the reader and amounts to checking the following relations 
 \begin{center}{\psfig{figure=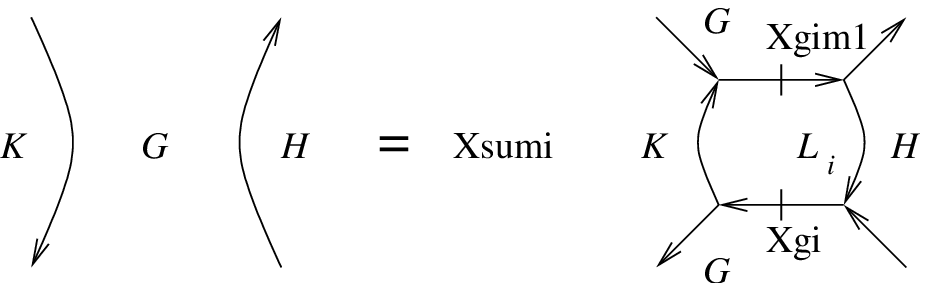}}\end{center}

 \begin{center}{\psfig{figure=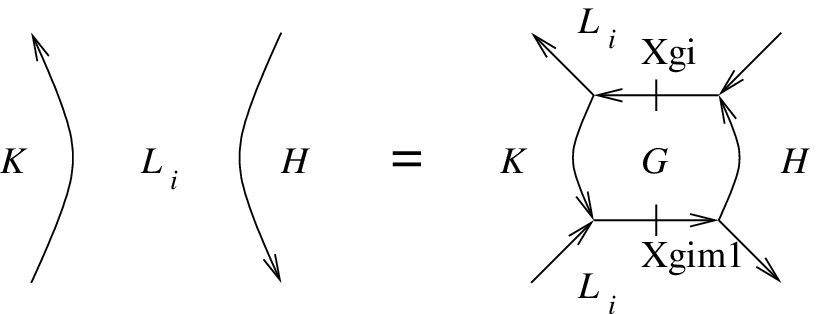}}\end{center}

\begin{center}{\psfig{figure=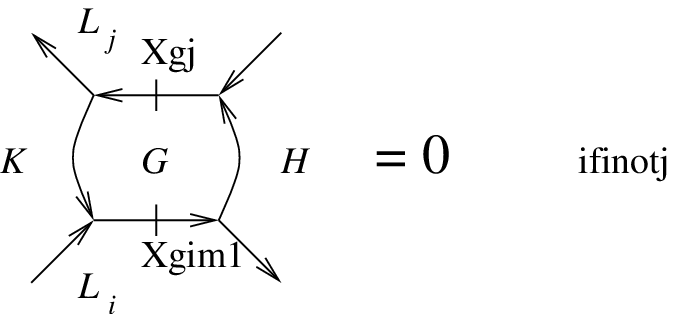}}\end{center}

%%%%%%%%%%%%%%%%%%%%%%%%
%%%%%%%%%%%%%%%%%%%%%%%
%%
%%  IND AND RES FOR SYMMETRIC GROUPS
%%
%%%%%%%%%%%%%%%%%%%%%%%%
%%%%%%%%%%%%%%%%%%%%%%%%%

\subsection{Induction and restriction between symmetric groups} \label{sub-indres} 

We now specialize the earlier construction to the case of the symmetric group $S_n$, viewed 
as the permutation group of $\{1,2,\dots, n\}$, and induction/restriction 
functors for inclusions $S_n \subset S_{n+1}$, where $S_n$ is identified with 
the stabilizer of $n+1$ in $S_{n+1}$. Notations for bimodules will be further 
simplified, so that, for instance, ${}_n (n+1)_{n-1}$ stands for $\Bbbk[S_{n+1}]$, 
viewed as a $(\Bbbk[S_n], \Bbbk[S_{n-1}])$-bimodule for the standard inclusions 
$S_n \subset S_{n+1}\supset S_{n-1}$, and ${}_n (n+1)_{n} (n+2)$ stands 
for $\Bbbk[S_{n+1}]\otimes_{\Bbbk[S_n]} \Bbbk[S_{n+2}]$, viewed 
as a $(\Bbbk[S_n], \Bbbk[S_{n+2}])$-bimodule. The regions of the strip $\R\times [0,1]$ 
are now labelled by nonnegative integers $n$. An upward oriented line separating 
regions labelled $n$ and $n+1$ 

\begin{center}{\psfig{figure=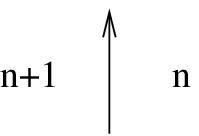}}\end{center}

\noindent 
denotes the identity endomorphism of the induction functor 
$$\Ind_n^{n+1} \  : \  \Bbbk[S_n]\dmod \lra \Bbbk[S_{n+1}]-mod .$$ 
This is the functor of tensoring with the bimodule $(n+1)_n$. 

A downward oriented line separating regions $n+1$ and $n$ 

\begin{center}{\psfig{figure=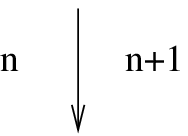}}\end{center}

\noindent 
denotes the identity endomorphism of the restriction functor 
$$\Res_{n+1}^{n} \  : \  \Bbbk[S_{n+1}]\dmod \lra \Bbbk[S_n]\dmod .$$
The bimodule for this functor is ${}_n(n+1)$. 

The four $U$-turns are given by the following bimodule maps: 

\begin{eqnarray}
\raisebox{-0.3cm}{\psfig{figure=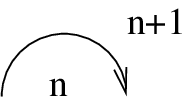}} & & 
\begin{array}{ccc} (n+1)_n(n+1) & \lra &  (n+1), \\    
 g\otimes h & \longmapsto & gh, \ g,h\in S_{n+1},\end{array} \\
\raisebox{-0.4cm}{\psfig{figure=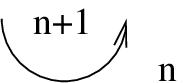}} & & 
(n) \lra {}_n(n+1)_n , \ \    g \longmapsto g, \  g\in S_n  \\
\raisebox{-0.2cm}{\psfig{figure=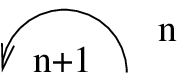}} & & 
p_n\ : \ {}_n(n+1)_n\lra (n) , \ \    p_n(g)=
\begin{cases} g & \text{if $g\in S_n$}, \\ 0 & \text{otherwise} 
\end{cases}  \\
\raisebox{-0.4cm}{\psfig{figure=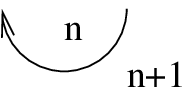}} & & 
q_n \ : \ (n+1)\lra (n+1)_n (n+1) , 
\end{eqnarray} 
where the bimodule map $q_n$ is determined by the condition 
$$ q_n(1) = \sum_{i=1}^{n+1} s_i s_{i+1}\dots s_n \otimes s_n \dots s_{i+1}s_i, 
\ \  \ s_i=(i,i+1).$$ 
Notice that $\{ s_n \dots s_2 s_1, s_n \dots s_3 s_2, \dots, s_n s_{n-1}, s_n, 1 \}$ 
are $n+1$ coset representatives of $S_n\subset S_{n+1}$, and 
\begin{eqnarray*} 
q_n(g) & = & \sum_{i=1}^{n+1} gs_i s_{i+1}\dots s_n \otimes s_n \dots s_{i+1}s_i \\
 & = &   \sum_{i=1}^{n+1} s_i s_{i+1}\dots s_n \otimes s_n \dots s_{i+1}s_i g, 
  \ \  g\in S_{n+1}.
\end{eqnarray*} 
Bimodule maps $p_n$ and $q_n$ are the second adjointness maps 
for the group inclusion $S_n\subset S_{n+1}$. For an arbitrary inclusion 
of finite groups $H\subset G$ these maps were described in the previous 
subsection. 

Denote by an upward-pointing crossing
 \raisebox{-0.4cm}{\psfig{figure=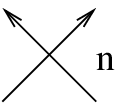}} the endomorphism 
of $(n+2)_n$ given by right multiplication by $s_{n+1}$, so that 
$g\longmapsto g s_{n+1}, g\in S_{n+2}$. 

Denote by  a downward-pointing crossing
 \raisebox{-0.4cm}{\psfig{figure=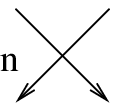}} the endomorphism 
of ${}_n(n+2)$ given by left multiplication by $s_{n+1}$, so that 
$g\longmapsto s_{n+1}g , g\in S_{n+2}$. 

Denote by  a right-pointing crossing
 \raisebox{-0.4cm}{\psfig{figure=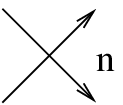}} the bimodule endomorphism 
$(n)_{n-1}(n) \ \lra  \  {}_n(n+1)_n$ that takes $g\otimes h$ for $g,h\in S_n$ to 
$g s_n h \in S_{n+1}$. 

Denote by  a left-pointing crossing
 \raisebox{-0.4cm}{\psfig{figure=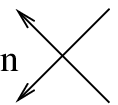}} the bimodule endomorphism 
${}_n(n+1)_n \  \lra \  (n)_{n-1}(n)$ that takes $g\in S_n\subset S_{n+1}$
to $0$ and $ g s_n h$ for $g,h\in S_n$ to $g\otimes h\in (n)_{n-1}(n)$. 
\vspace{0.2in} 

These four definitions-notations are compatible with the isotopies of diagrams 
in the plane strip--there are equalities of bimodule endomorphisms 

\begin{center}{\psfig{figure=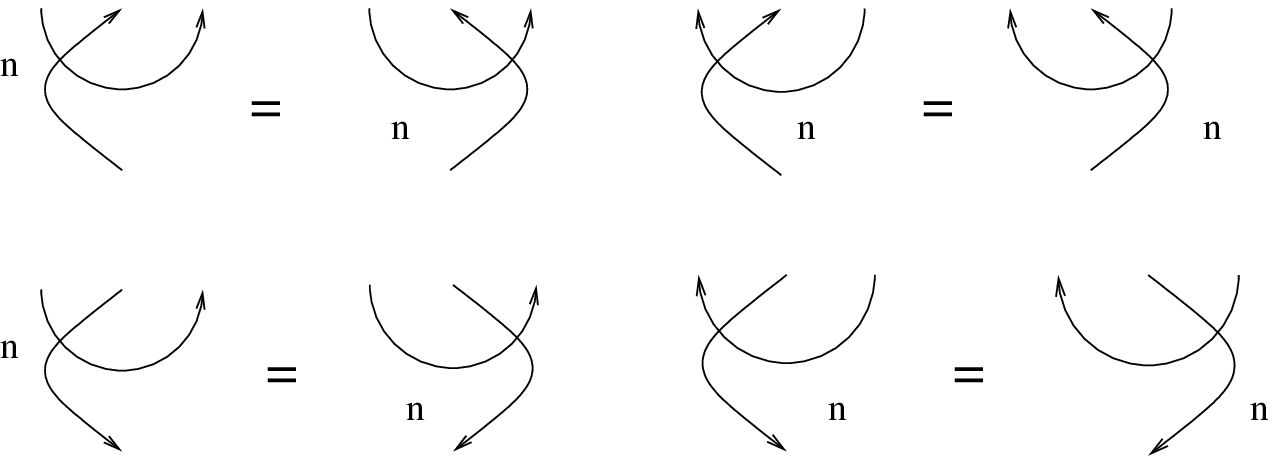}}\end{center}

They can be checked by direct computations. 

\begin{prop} \label{prop-6-again} The following relations hold for any $n\in \Z$. 
\begin{equation}\label{eq-mainrels4}
{\psfig{figure=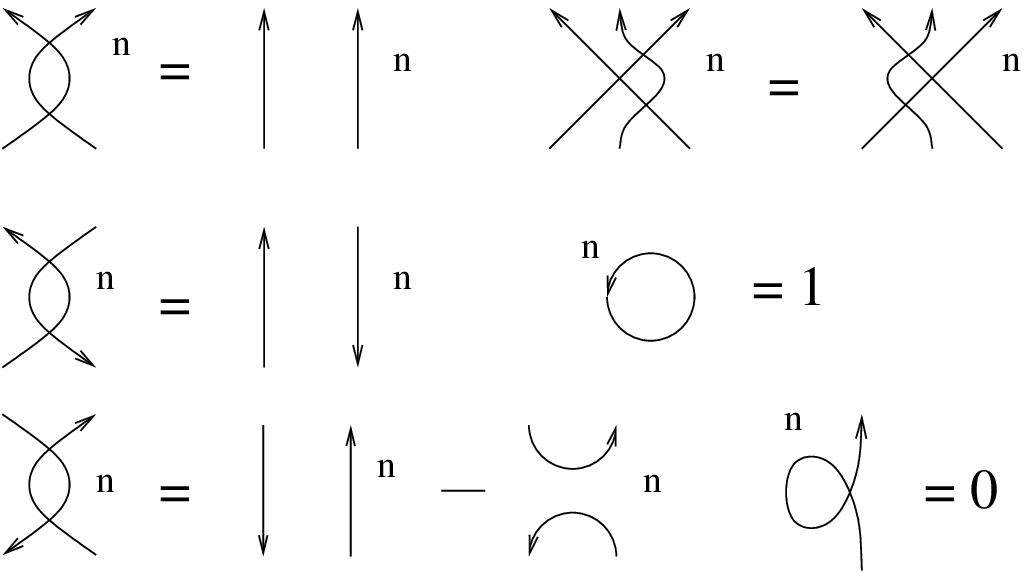}}\end{equation}
\end{prop} 

Here and further we follow the convention that a diagram equal $0$ if 
it has a region labelled by a negative number.  
Notice that the top two relations come from the relations in the symmetric groups:  
$s_{n+1}^2=1$ and $s_{n+1}s_{n+2}s_{n+1}=s_{n+2}s_{n+1}s_{n+2}$ 
and follow at once from the definition of the bimodule homomorphism 
associated to an upward-pointing crossing. The four remaining relations 
encode the bimodule decomposition 
$$ {}_n(n+1)_n \ \cong \  (n)_{n-1} (n) \oplus (n) $$ 
giving an isomorphism 
$$ \Res_{n+1}^n \circ \Ind_n^{n+1} \cong \Ind_{n-1}^n \circ 
\Res_n^{n-1} \oplus \Id$$ 
of endofunctors in the category $\Bbbk[S_n]\dmod$. This is a special case 
of the Mackey decomposition theorem which was given a diagrammatic 
interpretation in an earlier section. 
Proposition~\ref{prop-6-again} relations are identical to the ones in the definition of the 
category $\Cone,$ see Section~\ref{subsection-moves}.  

Let $\Sone$ be the category whose objects are compositions of 
induction and restriction functors between symmetric groups 
for standard embeddings $S_k \subset S_{k+1}$. The morphisms 
are natural transformations of functors (again, we work over a field 
$\Bbbk$ of characteristic $0$). 
The category $\Sone$ is the sum of categories $\Sone_k$ over $k\ge 0$; in 
the latter the first induction or restriction starts from $S_k$. 
For instance, $\Ind^{k+2}_{k+1} \circ \Ind^{k+1}_k \circ 
\Res_{k+1}^{k}\circ \Ind^{k+1}_{k}$ is an object of $\Sone_k$. 
Morphisms in $\Sone_k$ are natural transformation of functors
and can be identified with homomorphisms of associated bimodules. 

Thus, for each $k\ge 0$,  there is a functor $\mc{F}'_k : \Cone \lra \Sone_k$ 
that takes $Q_{\epsilon}$ to the corresponding composition of 
induction and restriction functors. For instance, 
$$\mc{F}'_k(Q_{++-+}) = \Ind^{k+2}_{k+1} \circ \Ind^{k+1}_k \circ 
\Res_{k+1}^{k}\circ \Ind^{k+1}_{k}. $$
If, for some $m$, the last $m$ terms of $\epsilon$ contain at least $k+1$ more 
minuses than pluses, then $\mc{F}'_k(Q_{\epsilon})=0$. On morphisms $\mc{F}'_k$ 
is defined as follows. It takes a diagram representing a morphism in $\Cone$, 
labels the rightmost region of the diagram by $k$, and views the diagram as 
a natural transformation between compositions of induction and restriction 
functors. 
The functor $\mc{F}'_k$ is not monoidal, since $\Sone_k$ does not 
have a monoidal structure matching that of $\Cone$. 

Let $\Stwo$, respectively $\Stwo_k$,  be the Karoubi envelope 
of $\Sone$, respectively $\Sone_k$. Functor $\mc{F}'_k$ induces a 
functor on Karoubi envelopes $ \mc{F}_k \ : \ \Ctwo \lra \Stwo_k$.
We summarize relevant categories and functors below. 

\begin{equation*} 
\Stwo = \mathrm{Kar}(\Sone), \ \  \Stwo_k = \mathrm{Kar}(\Sone_k), \ \ 
  \mc{F}_k' : \Cone \lra \Sone_k, \ \  \mc{F}_k : \Ctwo \lra \Stwo_k . 
\end{equation*} 

$$ 
\begin{CD} 
\Sone  @=   \oplusop{k\ge 0} \Sone_k  \\
 @V{\mathrm{Kar}}VV       @VV{\mathrm{Kar}}V  \\
 \Stwo    @=   \oplusop{k\ge 0} \Stwo_k 
\end{CD}   \ \ \ \ \ \ \ \ \ \ \ \ \ \ \ \ \ \ 
\begin{CD} 
\Cone  @>{\mc{F}_k'}>> \Sone_k \\
 @V{\mathrm{Kar}}VV       @VV{\mathrm{Kar}}V  \\
 \Ctwo  @>{\mc{F}_k}>> \Stwo_k
\end{CD} 
$$

Functor $\mc{F}_k$ induces a homomorphism of abelian groups 
\begin{equation} 
[\mc{F}_k] \ : \ K_0(\Ctwo) \lra K_0(\Stwo_k).  
\end{equation} 
Notice that $K_0(\Ctwo)$ is a ring, while $K_0(\Stwo_k)$ is only an 
abelian group. 

An object  of $\Stwo_k$ is a direct summand of a finite sum of composition 
of induction and restriction functors that start with the category of $\Bbbk[S_k]$-modules, 
thus it takes any finite-dimensional  $\Bbbk[S_k]$-module  to a module over  
$\oplusop{m\ge 0} \Bbbk[S_m]$. Descending to Grothendieck groups, we obtain 
a homomorphism 
$$ \theta_k \ : \ K_0(\Stwo_k) \lra \Hom_{\Z}(K_0(\Bbbk[S_k]),
 \oplusop{m\ge 0} K_0(\Bbbk[S_m])) .$$ 

From here until the end of this paper we assume that ${\mathrm char}(\Bbbk)=0$.
Consider the composite homomorphism 
$$ \theta_k [\mc{F}_k] \ : \ K_0(\Ctwo) \lra \Hom_{\Z}(K_0(\Bbbk[S_k]),
 \oplusop{m\ge 0} K_0(\Bbbk[S_m])) .$$ 
This homomorphism takes $[Q_{+,\mu}]$ to a map that assigns 
to $[M]\in K_0(\Bbbk[S_k]),$ for a $\Bbbk[S_k]$-module $M$, the 
symbol  
$[\mathrm{Ind}_{S_{|\mu|}\times S_k}^{S_{|\mu|+k}}(L_{\mu}\otimes M)] $
of induced module over $\Bbbk[S_{|\mu|+k}]$. In other words, tensor $M$ with 
$L_{\mu}$, producing a module over $\Bbbk[S_{|\mu|}]\times S_k$, induce 
to $\Bbbk[S_{|\mu|+k}]$, then pass to the Grothendieck group. 

Likewise, $\theta_k [\mc{F}_k]$ takes $[Q_{-,\lambda}]$ to the zero map if 
$|\lambda|> k$ and, if $k \ge |\lambda|,$  
to the map that assigns to $[M]$ as above the symbol of the module 
$$ \Hom_{\Bbbk[S_{|\lambda|}]}(L_{\lambda}, M) \in \Bbbk[S_{k-|\lambda|}]\dmod$$ 
In other words, restrict $M$ to being a module over the group algebra of  
$S_{|\lambda|}\times S_{k-|\lambda|} \subset S_k$, and form homs from the simple 
module $L_{\lambda}$ over $S_{|\lambda|}$. The result is a representation of 
the symmetric group $S_{k-|\lambda|}.$ 

Now consider composition 
$$\theta_k [\mc{F}_k]\gamma \ : \ H_{\Z} \lra 
\Hom_{\Z}(K_0(\Bbbk[S_k]), \oplusop{m\ge 0} K_0(\Bbbk[S_m])).$$ 

We claim that the sum of these maps, over all $k\ge 0$, is injective. Let 
$$y = \sum_{\lambda,\mu} 
y_{\lambda, \mu} b_{\mu} a_{\lambda}, \ \ y_{\lambda,\mu}\in \Z$$  
be an arbitrary nonzero element of $H_{\Z}$. We have 
 $$ \gamma(y) = \sum_{\lambda,\mu} y_{\lambda,\mu} [Q_{+,\mu^{\ast}}]
[Q_{-,\lambda}].$$ 
When $|\lambda|=k$, the map 
$$\theta_k [\mc{F}_k]\gamma(a_{\lambda}) = \theta_k [\mc{F}_k]([Q_{-,\lambda}])$$ 
takes $[L_{\nu}]$ to $0$ if $|\nu|=k$ and $\nu\not= \lambda$. 
The same map takes $[L_{\lambda}]$ to $[L_{\emptyset}]$, the symbol of the simple 
module over $\Bbbk[S_0] = \Bbbk$. 

Choose $k$ such that $y_{\lambda,\mu}\not= 0$ for some $\lambda$ with 
$|\lambda|=k$ and some $\mu$, while $y_{\lambda,\mu}=0$ for all $\mu$ 
whenever $|\lambda|<k$. 
Also choose $\nu$ with $|\nu|=k$ and $y_{\nu,\mu}\not= 0$ for some 
$\mu$.  
The map $\theta_k [\mc{F}_k]\gamma(y)$ takes 
$[L_{\nu}]$ to 
$$\sum_{\mu} y_{\nu,\mu} [L_{\mu^{\ast}}] \not= 0.$$  
Therefore, $\theta_k [\mc{F}_k]\gamma(y)$ is a nonzero map, and $\gamma(y)\not= 0$. 
This concludes the proof that $\gamma$ is injective (Theorem~\ref{thm-inj}).

%%%%%%%%%%%%%%%%%%%%%%%%%
%%%%%%%%%%%%%%%%%%%%%%%%%
%%
%%  SIZE OF C 
%%
%%%%%%%%%%%%%%%%%%%%%%%%%
%%%%%%%%%%%%%%%%%%%%%%%%%

\section{The size of morphism spaces in $\Cone$} \label{sec-sizeofC}

In this section we will prove Propositions~\ref{prop-psi-iso} and 
\ref{prop-psim-iso}. 
Consider the  right curl with the rightmost region labelled $n$ (also recall 
the shorthand of denoting this curl by a dot). This curl can be realized as 
the composition of a cup with a crossing with a cap.  

\begin{center}{\psfig{figure=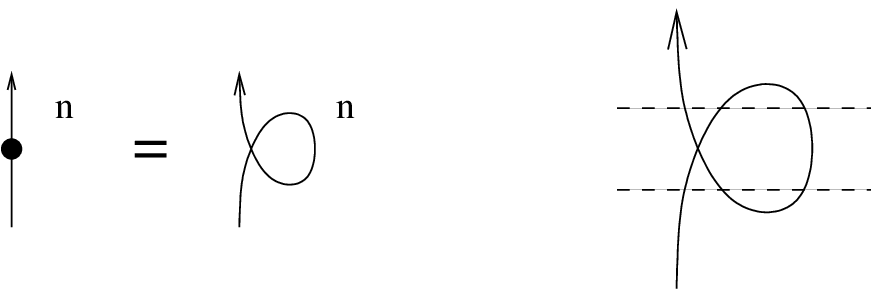}}\end{center}

The corresponding endomorphism of the bimodule $(n+1)_n$ takes $1$ to 
$$J_n :=\sum_{i=1}^n s_i \dots s_{n-1}s_n s_{n-1} \dots s_i=(1,n+1)+(2,n+1)+\dots 
+(n,n+1). $$ 
This endomorphism of $(n+1)_n$ is the right multiplication by $J_n$: 
$$ g \longmapsto g J_n, \ \  g\in S_{n+1}.$$
Notice that $J_n$ is the Jucys-Murphy element, ubiquitous in the representation 
theory of the symmetric group. Our diagrammatics realizes it via the 
right  curl and interprets the endomorphism of multiplication by $J_n$ 
as the composition of three natural transformations, two of which (cup and cap) 
come from biadjointness of the induction and restriction functors. 

Commutativity of Jucys-Murphy elements now acquires a graphical interpretation 
as isotopies of curls (or dots) on upward-oriented disjoint strands past each other.  

\begin{center}{\psfig{figure=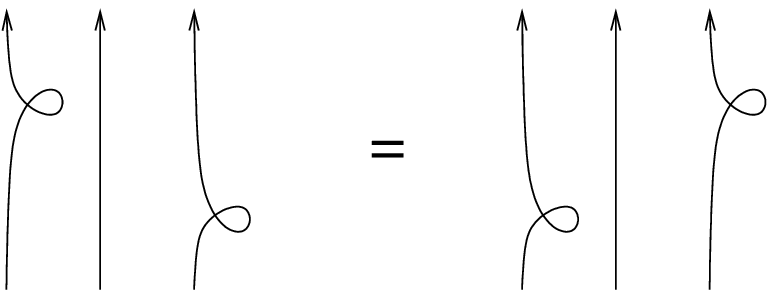}}\end{center}

For each $n \ge 0$ the functor $\mc{F}'_n$, applied to $\mathbf{1}$ and 
its endomorphisms, produces  a homomorphism from 
$\End_{\mc{H}'}(\mathbf{1})$ to $Z(\Bbbk[S_n])$, the center of the group ring of the 
symmetric group. Composing with the homomorphism $\psi_0$, we obtain 
the homomorphism 
$$\psi_{0,n}: \Pi \cong \Bbbk [c_0, c_1, \dots ] \lra Z(\Bbbk [S_n]) $$ 
that takes $c_k$ to
\begin{eqnarray} 
c_{k,n} & = & \sum_{i=1}^n s_i \dots s_{n-1} J_{n-1}^k s_{n-1} \dots s_i \\
    & = & \sum_{i=1}^n (i,i+1,\dots, n)\big(  (1,n)+\dots + (n-1,n)\big)^k 
 (n, n-1, \dots, i).  
\end{eqnarray} 
We'd like to show that the union of $\psi_{0,n}$ over all $n$, 
$$ \Pi \lra \oplusop{n\ge 0}Z(\Bbbk [S_n]), $$
is injective. For the first few values of $k$ we have 
\begin{eqnarray*} 
 c_{0,n} & = & n , \\
 c_{1,n} & = &  2 \sum_{1\le i < j \le n} (i,j) , \\
 c_{2,n} & = &  3 \sum (i_1, i_2, i_3) +  n(n-1) , \\
 \end{eqnarray*}  
where the sum in the bottom line is over all 3-cycles. In general, 
\begin{equation} \label{eq-cnk}
c_{k,n} = (k+1)\sum (i_1, \dots, i_{k+1}) + \mathrm{l.o.t.}, 
\end{equation} 
where the sum is over all $k+1$-cycles. For $k>0$, the lower order terms 
is a sum over permutations of disturbance less than $k+1$, 
where we define the disturbance of a permutation $\sigma\in S_n$ 
as the number of elements moved by $\sigma$, 
$$dist(\sigma) = | \{ i | 1\le i\le n, \sigma(i) \not= i\} |.$$ 
 Notice that 
$c_{k,n}$ is the sum of conjugates of $J_{n-1}^k$ and contains terms of 
disturbance at most $k+1$, for $k>0$. 

Since ${\mathrm char}(\Bbbk)=0$,  the 
coefficient of the sum of $(k+1)$-cycles in (\ref{eq-cnk}) is nonzero. 
Consider an increasing filtration 
$$ Z_0=\Bbbk \subset Z_2 \subset Z_3 \subset \dots \subset Z_n = \Z(\Bbbk[S_n])$$ 
by having $Z_i$ be the span over conjugacy classes that consist of permutations 
of disturbance at most $i.$ 
We make $\Bbbk[c_0, c_1,\dots ]$ graded by $\deg(c_0)=0$ and $\deg(c_k)=k+1$ 
for $k>0$, and then consider associated increasing filtration on $\Bbbk[c_0, c_1,\dots ]$, 
where the $i$-th filtered terms is the sum of graded terms of degree at most $i$. 
The homomorphism $\psi_{0,n}$ respects the two filtrations.

To show asymptotic faithfulness of $\psi_{0,n}$ as $n \lra \infty$ we examine 
induced homomorphism of adjoint graded rings. Assume that 
there is a universal relation in $\End_{\mc{H}'}(\mathbf{1})$ 
$$ \sum_{I}  a_I c_I =0 $$ 
which holds for all sufficiently large $n$. Here $I=(a_0,a_1,\dots, a_i)$ 
is a finite sequence of non-negative integers, $a_I\in \Bbbk$, and 
$c_I = c_0^{a_0}c_1^{a_1} \dots c_i^{a_i}$ is the corresponding monomial 
in $\Bbbk[c_0,c_1, \dots ]$. The sum is over finitely many sequences $I$. 
This implies relations 
$$ \sum_{I}  a_I c_{I,n} =0 $$ 
for all $n$, where $c_{I,n} = c_{0,n}^{a_0}c_{1,n}^{a_1} \dots c_{i,n}^{a_i}$. 
Choose any term from the sum corresponding to $I=(a_0,a_1,\dots, a_i)$ with 
the highest possible degree monomial $c_I$. The only terms in the 
sum $\sum_{I}  a_I c_{I,n}$ that contribute to the conjugacy class of type 
$(a_1+1, \dots, a_i+1)$ can come from sequences 
$I'=(a, a_1, \dots, a_i)$ that differ from $I$ only in the first term. 
Since, for such $I'$, 
$$ c_{I',n}= n^a c_{1,n}^{a_1}\dots c_{i,n}^{a_i} $$ 
the relation $\sum_{I}  a_I c_{I,n} =0 $ implies 
$\sum_a  a_{I'} n^a =0$, which, in turn, leads to $a_{I'}=0$ for all $I'$ as 
above, since the only polynomial with infinitely many positive integers $n$ as roots is 
the zero polynomial. This contradiction implies asymptotic faithfulness of 
$\psi_{0,n}$. 
Therefore, $\psi_0$ is injective, which concludes the proof of 
Proposition~\ref{prop-psi-iso} when ${\mathrm char}\Bbbk=0$. 
Take $\Bbbk = \Q$, then pass to the subring $\Z$ of $\Q$. Our arguments imply 
that $\psi_0$ is an isomorphism over $\Z$, and, hence, over any commutative 
ring $\Bbbk$, including any field. 

\vspace{0.1in} 

We next prove injectivity of $\psi_m$  (Proposition~\ref{prop-psim-iso}),  
see equation (\ref{eq-psim}), again first in the characteristic zero case.  
The algebra $DH_m\otimes \Pi$ has a basis of elements 
$x_1^{a_1}\dots x_m^{a_m} \cdot \sigma \cdot c_0^{b_0} c_1^{b_1}\dots c_k^{b_k}$
over permutations $\sigma\in S_m$ and $a_i, b_j\in \Z_+,$ see discussion 
before Proposition~\ref{prop-psim-iso}. 

Applying functor $\mc{F}_n'$ to $Q_{+^m}$ and its endomorphism ring gives 
us a homomorphism 
$$ \End_{\Cone}(Q_{+^m}) \lra \End((n+m)_n) $$ 
to the endomorphism ring of the $(\Bbbk[S_{n+m}], \Bbbk[S_n])$-bimodule 
$\Bbbk[S_{n+m}]$, which we also denote $(n+m)_n$. 
The composite homomorphism 
$$ \psi_{m,n} \ : \ DH_m \otimes \Pi \lra  \End((n+m)_n)  $$ 
takes elements of $\ DH_m \otimes \Pi$ to endomorphisms given by right 
multiplication by suitable elements of $\Bbbk[S_{n+m}]$. 
Namely, $\psi_{m,n}(\sigma)$, for a permutation $\sigma\in S_m \subset DH_m$,  
is the right multiplication by $\underline{\sigma}$, where we define the latter 
as the translate of $\sigma$ by $n$: 
$$ \underline{\sigma}(i+n) = \sigma(i)+n, \ \ \  1\le i \le m, \ \ 
   \underline{\sigma}(i) = i, \ \ \  1\le i \le n.$$ 
$\psi_{m,n}(x_i)$, where  $x_i$ the the diagram of $m$ vertical lines 
with the dot (right  curl) on the $i$-th strand from the left, is the endomorphism 
of right multiplication by $J_{n+m-i}$, and $\psi_{m,n}(c_k)$ is the right 
multiplication by $c_{k,n}$.  The map $\psi_{m,n}$ is described by 
the corresponding homomorphism $\psi'_{m,n}$ from $DH_m\otimes \Pi$ 
to the opposite of the group algebra $\Bbbk[S_{n+m}]^{op}\supset \End((n+m)_n)$, 
with $\psi'_{m,n}(x_i) = J_{n+m-i}$, $\psi'_{m,n}(c_k) = c_{k,n}$, etc.  
We need to take the opposite algebra since the ring of endomorphism of a ring $A$ 
viewed as a left $A$-module is naturally isomorphic to the opposite of $A$: 
 $ \End_A ({}_A A, \ {}_A A) \cong A^{op}.$

Define $m$-disturbance of a permutation $\sigma\in S_{n+m}$ as 
the number of integers between $1$ and $n$ that are moved by $\sigma$: 
$$ dist_m(\sigma) = |\{ i | 1\le i \le n , \sigma(i) \not= i \} |.$$ 
Notice that $dist_m(\sigma\tau) \le dist_m(\sigma) + dist_m(\tau)$. 
On the group algebra $\Bbbk[S_{n+m}]$ we can introduce an increasing 
filtration 
$$\Bbbk[S_m] = Z_0^m \subset Z_1^m \subset \dots \subset Z^m_n = \Bbbk[S_{n+m}]$$ 
where $Z_k^m$ is spanned by permutations of disturbance at most $k$. 

We turn $DH_m\otimes\Pi$ into a filtered algebra 
by setting $\deg(c_0)=0, $ 
$\deg(c_k) = k+1$ if $k>0$, $\deg(x_i) = 1$ and $\deg(\sigma)=0$ 
and then making the $k$-th term in the increasing filtration be spanned 
by the basis elements of total degree at most $k$. 

Homomorphism $\psi'_{m,n}$ is that of filtered algebras, and we can 
pass to the homomorphism of associated graded algebras. 
To show asymptotic faithfulness of $\psi'_{m,n}$ we fix $m$ and 
will be taking $n$ large compared to $m$. Assume that 
there exists a relation 
\begin{equation}\label{eq-rel0}
\sum d_{\sigma,{\bf a},{\bf b}}  \ \ 
x_1^{a_1}\dots x_m^{a_m} \cdot \sigma \cdot c_0^{b_0} c_1^{b_1}\dots c_r^{b_r} 
=0
\end{equation} 
in $\End_{\Cone}(Q_{+^m})$ for some $d_{\sigma, {\bf a}, {\bf b}}\in \Bbbk\setminus 
\{0\}$, with ${\bf a}=(a_1,\dots, a_m),$ ${\bf b}= (b_1, \dots, b_r)$, the 
sum over finitely many triples $( \sigma,{\bf a},{\bf b})$. 
Let 
$$x(\sigma, {\bf a}, {\bf b}) = x_1^{a_1}\dots x_m^{a_m} \cdot
\sigma  \cdot c_0^{b_0} c_1^{b_1}\dots c_r^{b_r} $$
denote the elements of our basis of $DH_m\otimes \Pi$.  
The element  
$ \psi'_{m,n}( x(\sigma, {\bf a}, {\bf b}))\in \Bbbk[S_{n+m}]$ belongs 
to the $k$-th term of the filtration of the latter, where 
$$ k = a_1 + \dots + a_m + 2 b_1 + 3b_2 + \dots +(r+1)b_r,$$
but not to the $(k-1)$-st term. 
Among terms $ x(\sigma, {\bf a}, {\bf b})$ in the sum select only those with the 
maximal possible $k$ (denote such $k$  by $k_0$). 
It's enough to show that, as we sum over only these terms, 
the image of 
$ \sum \psi'_{m,n} ( x(\sigma, {\bf a}, {\bf b})) $ in the associated graded ring 
of $\Bbbk[S_{n+m}]$ relative to the above filtration is nonzero. In other 
words, we need to show that coefficients of permutations of disturbance 
$k_0$ are not all zero for some sufficiently large $n$ in the expansion 
of $\psi_{m,n}$ applied to the LHS of (\ref{eq-rel0}). This is obtained 
by looking at the structure of these permutations. They are disjoint unions 
of cycles, with some of the cycles containing one or more elements of  
the set $P=\{ n+1, \dots, n+m\}$.  The relative positions of elements 
of $P$ in the cycles, lengths of portions of the cycles between elements 
of $P$, together with  the number of cycles of each length 
without elements of $P$ uniquely determine $\sigma$, $a_1, \dots , a_m$ 
and $c_1, \dots, c_r$ that can contribute to the coefficient of such permutation. 
Coefficients at different powers  of $c_0$ are 
taken care in the same way as in the $m=0$ case. 
Linear independence of our spanning set of $\End_{\Cone}(Q_{+^m})$ 
and Proposition~\ref{prop-psim-iso} follow when ${\mathrm char}\Bbbk=0$. 
Same argument as in the $m=0$ case then implies that $\psi_m$ is 
an isomorphism over any commutative ring $\Bbbk$.  

\vspace{0.14in} 

The formula (\ref{eq-cnk}) implies that the natural homomorphism 
$\End_{\Cone}({\mathbf 1})\lra Z(\Bbbk[S_n])$ 
from the endomorphisms of the identity object of $\Cone$ to the center 
of the group algebra is surjective when the field $\Bbbk$ has characteristic 
$0$. Combining with the  result of Cherednik~\cite{Cher2} and Olshanskii~\cite{OL}, 
\cite[Theorem 3.2.6]{CST}
that the centralizer algebra of $\Bbbk[S_n]$ in $\Bbbk[S_{n+m}]$ 
is generated by $DH_m$ and the center of $\Bbbk[S_n]$, we obtain that 
the homomorphism $DH_m\otimes \Pi \lra \End((n+m)_n)$ introduced 
above is surjective when ${\mathrm char} \Bbbk = 0$.

%%%%%%%%%%%%%%%%%%%%%%%%%
%%%%%%%%%%%%%%%%%%%%%%%%%
%%
%%  REMARKS 
%%
%%%%%%%%%%%%%%%%%%%%%%%%%
%%%%%%%%%%%%%%%%%%%%%%%%%

\section{Remarks on the Grothendieck ring of $\Ctwo$} \label{sec-sizingup}

%%%%%%%%%%%%%%%%%%%
%%%%%%%%%%%%%%%%%%%
%
%  IDEMPOTENTED RINGS 
%
%%%%%%%%%%%%%%%%%%%
%%%%%%%%%%%%%%%%%%%

\subsection{Idempotented rings from $\Cone$} 

For a sequence $\epsilon$ of pluses and minuses denote by $\langle \epsilon \rangle$ 
the difference between the number of pluses and minuses in $\epsilon$ (the weight of 
$\epsilon$).  Then 
$\Hom_{\Cone}(Q_{\epsilon}, Q_{\epsilon'}) = 0$ if and only if 
$\langle \epsilon \rangle \not= \langle \epsilon' \rangle . $ 
The "if" part of this observation 
implies that categories $\Cone$ and $\Ctwo$, viewed as additive categories, decompose 
into the direct sum of subcategories 
$$ \Cone = \oplusop{\ell\in\Z} \Cone_{\ell}, \ \ \ \ \ 
     \Ctwo = \oplusop{\ell\in\Z} \Ctwo_{\ell}, $$
where $\Cone_{\ell}$ is a full subcategory of $\Cone$ which contains objects 
$Q_{\epsilon}$ over all sequences of weight $\ell$, 
and $\Ctwo_{\ell}$ is the Karoubi envelope of $\Cone_{\ell}$. 

This direct sum decomposition induces a grading on the Grothendieck 
ring 
$$K_0(\Ctwo)= \oplusop{\ell \in \Z} K_0(\Ctwo_{\ell})$$ 
The Heisenberg algebra $H$ and its integral form $H_{\Z}$ are graded by 
$\deg(a_n) = n = - \deg(b_n)$, and the homomorphism 
$\gamma: H_{\Z} \lra K_0(\Ctwo)$ is that of graded rings.

We can redefine the Grothendieck groups $K_0(\Ctwo)$ and $K_0(\Ctwo_{\ell})$ via 
idempontented rings. 
For a sequence $\epsilon$ let 
$$ \End(\epsilon) \ := \ \End_{\Ctwo}(Q_{\epsilon})$$ 
denote the endomorphism algebra of $Q_{\epsilon}$. Likewise, denote 
$$ \Hom(\epsilon, \epsilon') \ := \ \Hom_{\Ctwo}(Q_{\epsilon}, Q_{\epsilon'}).$$  
To the category $\Ctwo$ we can assign the idempotented ring of all homomorphisms 
between various tensor products of generating objects $Q_+$ and $Q_-$:  
$$ \underline{R} \ := \oplusop{\epsilon,\epsilon'} \Hom(\epsilon, \epsilon'),$$ 
the sum over all sequences $\epsilon, \epsilon'$. Ring $\underline{R}$ is nonunital, but has 
a family of idempotents $1_{\epsilon} =1 \in \End(\epsilon).$
 
Right projective $\underline{R}$-modules 
$1_{\epsilon}\underline{R}$ correspond to objects $Q_{\epsilon}$, 
in the sense that 
\begin{equation} \label{eq-twoeq}
\Hom_{\underline{R}}(1_{\epsilon}\underline{R}, 1_{\epsilon'}\underline{R}) = 
\Hom(\epsilon, \epsilon') = 
\Hom_{\Cone}(Q_{\epsilon}, Q_{\epsilon'}),
\end{equation} 
and the Grothendieck group of finitely-generated projective right $\underline{R}$-modules 
is canonically isomorphic to the Grothendieck group of $\Ctwo$: 
 $$ K_0(\underline{R}) \cong K_0(\Ctwo).$$ 
This isomorphism takes $[1_{\epsilon}\underline{R}]$ to $[Q_{\epsilon}]$. 
Usually we use $K_0(A)$ to denote the Grothendieck group of finitely-generated 
projective left, not right, $A$-modules. Here, because of (\ref{eq-twoeq}), we 
use right $\underline{R}$-modules in the definition of $K_0(\underline{R})$. 
Alternatively, we could use 
left $\underline{R}^{op}$-modules, or even left $\underline{R}$-modules after 
fixing an isomorphism $\underline{R}\cong \underline{R}^{op}$ 
(involution $\xi_2$ induces one such isomorphism). We have 

$$ \underline{R} = \oplusop{\ell \in \Z} \underline{R}_{\ell} , \ \ \ \ \ 
  \underline{R}_{\ell} \ :=  \ \oplusop{\langle \epsilon\rangle =
\langle\epsilon'\rangle = \ell } \Hom(\epsilon, \epsilon').$$ 

Assume from now on that $\ell\ge 0$ (the other case can be treated similarly, 
or by applying symmetry $\xi_3$ to reverse $\ell$).
Given a sequence $\epsilon$ with $n+\ell $ pluses and $n$ minuses, 
the object $Q_{\epsilon}$ of $\Ctwo_{\ell}$ decomposes into the direct sum of 
objects $Q_{+^{k+\ell} -^{k}}$ with $0\le k \le n$ with some multiplicities. 
Hence, $\underline{R}_{\ell}$ is Morita equivalent to the idempotented ring 
$$ R_{\ell} \ = \ \oplusoop{k,k'=0}{\infty} \Hom(+^{k+\ell} -^{k}, 
 +^{k'+\ell} -^{k'}), $$ 
and the inclusion $R_{\ell}\subset \underline{R}_{\ell}$ induces an isomorphism 
of Grothendieck groups $K_0(R_{\ell}) \cong K_0(\underline{R}_{\ell}).$
Let 
$$ R_{\ell,m} \ = \ \oplusoop{k,k'=0}{m} \Hom(+^{k+\ell} -^{k}, 
 +^{k'+\ell} -^{k'}).$$
The ring $R_{\ell}$ is the union of rings in the increasing chain  
$R_{\ell,0}\subset R_{\ell,1}\subset\dots$. Formation of Grothendieck group 
commutes with direct limits,  implying that $K_0(R_{\ell})$ is the direct limit of 
$K_0(R_{\ell, m})$ as $m$ goes to infinity. Thus, there is an isomorphism 
$$ K_0(\Ctwo_{\ell}) \cong \lim_{m \to \infty} K_0(R_{\ell, m}).$$ 

For each $k$ between $0$ and $m$ natural inclusion of rings 
$\End(+^{k+\ell}-^k) \subset R_{\ell, m}$ induces a homomorphism 
of groups $K_0(\End(+^{k+\ell}-^k))\lra K_0(R_{\ell, m})$. 
Conjecture~\ref{conj-iso} would follow from the following two conjectures: 

\vspace{0.06in} 

\noindent 
{\it Conjecture 1.1.} The standard inclusion 
$\Bbbk[S_n\times S_m]\subset \End(+^n -^m)$ 
induces an isomorphism of Grothendieck groups of these two rings. 

\hspace{0.04in} 

\noindent 
{\it Conjecture 1.2.} Ring inclusion 
 $$ \oplusoop{k=0}{m}  \End(+^{k+\ell}-^k) \subset R_{\ell,m}$$ 
induces an isomorphism on Grothendieck groups. 

\vspace{0.06in}  

We don't know how to prove either statement, but will present now some 
weak evidence in favor of Conjecture 1.1.

%%%%%%%%%%%%%%%%%%%%%%%%%%%%%%%
%%%%%%%%%%%%%%%%%%%%%%%%%%%%%%%%
%%
%% GROTH GROUPS degenerate AHA 
%%
%%%%%%%%%%%%%%%%%%%%%%%%%%%%%%%%
%%%%%%%%%%%%%%%%%%%%%%%%%%%%%%%

\subsection{Grothendieck group of degenerate affine Hecke algebra}

Here we prove Conjecture 1.1 in the case $m=0$ ($n=0$ case follows by 
symmetry).  By Proposition~\ref{prop-psim-iso}, the endomorphism ring of the 
object $Q_{+^n}$ of $\Ctwo$ is isomorphic 
to the tensor product of the degenerate affine Hecke algebra $DH_n$ and 
the polynomial algebra $\Pi$: 
$$ \End(+^n) \cong DH_n \otimes   \Pi.$$ 
The inclusion $\Bbbk[S_n]\hookrightarrow DH_n$ is split, via 
the homomorphism $\tau_n: DH_n \lra \Bbbk[S_n]$ which takes 
generators $s_i$ of $DH_i$ to transpositions $(i,i+1)$ and 
generators $x_i$ to Jucys-Murphy elements. The split inclusion 
induces a split short exact sequence of two rings and an ideal 
$$ 0 \lra \mathrm{ker}(\tau_n) \lra DH_n \lra \Bbbk[S_n] \lra 0, $$
which, in turn, induces a split short exact sequence of $K_0$-groups 
$$ 0 \lra K_0(\mathrm{ker}(\tau_n)) \lra K_0(DH_n) \lra K_0(\Bbbk[S_n]) \lra 0,$$
see~\cite{Ros, Wei}.  Introduce an increasing filtration 

$$0 = Z_{-1}DH_n \subset Z_0 DH_n \subset Z_1 DH_n \subset \dots$$ 
on $DH_n$, where $Z_k DH_n$ is spanned by elements of the 
form $x_1^{a_1}\dots x_n^{a_n} \sigma$ over all $\sigma\in S_n$ and 
$a_1+\dots + a_n \le k$. Then $Z_k DH_n \times Z_m DH_n \subset Z_{k+m} DH_n$ 
and $Z_0 DH_n = \Bbbk [S_n]$. Let $B= \mathrm{gr} DH_n$ with respect to 
this filtration. $B$ is a graded algebra isomorphic to the cross-product 
of the polynomial algebra in $n$ generators with the group algebra of the 
symmetric group, $B\cong \Bbbk[x_1, \dots, x_n] \ast \Bbbk [S_n]$. 

Algebra $B$ is Koszul, with the Koszul dual algebra isomorphic to 
the cross-product of the exterior algebra on $n$ generators with the 
group algebra of the symmetric group (recall that $\mathrm{char}(\Bbbk)=0$).
Hence, $B$ has finite Tor dimension and, in particular, $Z_0 DH_n = \Bbbk[S_n]$ 
has Tor dimension $n$ as a right $B$-module. $B$ has Tor dimension $0$ 
as a right $Z_0 DH_n$-module, since $\Bbbk[S_n]$ is semisimple.  
We are in a position to invoke Quillen's theorem~\cite[Theorem 7, page 112]{Qui}, 
\cite{Block}.  

\vspace{0.06in} 

{\bf Theorem.} {\emph Let $A$ be a ring equipped with an increasing filtration 
$Z_k A$, and such that $Z_0 A$ is regular. Suppose that $B= \mathrm{gr}(A)$ 
has finite Tor dimension as a right $Z_0 A$ module and that $Z_0 A$ has 
finite Tor dimension as a right $B$ module. Then the inclusion $Z_0A \subset A$ 
induces an isomorphism $K_i (Z_0 A) \cong K_i(A)$.}  

\vspace{0.06in} 

In our case $A= DH_n$. The regularity of $Z_0 DH_n$ is obvious due to it 
being semisimple (a regular ring is a noetherian ring such that every left 
module has finite projective dimension). Applying the theorem in $i=0$ case 
we obtain 

\begin{prop} The inclusion $\Bbbk[S_n] \subset DH_n$ and the surjection 
$DH_n \longrightarrow \Bbbk[S_n]$ induce mutually-inverse  
isomorphisms $K_0(\Bbbk[S_n]) \cong K_0 (DH_n)$. 
\end{prop}   

Same argument shows that inclusion 
$$ \Bbbk[S_n] \subset DH_n \otimes \Bbbk[c_0, \dots, c_r]$$ 
induces an isomorphism of Grothendieck groups 
$$ K_0( \Bbbk[S_n]) \cong K_0(DH_n \otimes \Bbbk[c_0, \dots, c_r]).$$ 
Formation of Grothendieck groups commutes with taking direct limit of 
rings~\cite[Section 1.2]{Ros}, and $\Pi$ is the limit of $\Bbbk[c_0, \dots, c_r]$ 
as $r\to \infty$. We obtain an isomorphism 
$$ K_0( \Bbbk[S_n]) \cong K_0(\End(+^n))$$ 
proving Conjecture 1.1 when $m=0$.  $\square$ 

\vspace{0.06in} 

Induction and restriction functors for inclusions 
\begin{eqnarray*} 
  \Bbbk[S_n] \otimes \Bbbk[S_m] & \subset & \Bbbk[S_{n+m}]  \\
   DH_n \otimes DH_m  & \subset & DH_{n+m} 
\end{eqnarray*} 
induce ``multiplication'' and ``comultiplication'' maps on Grothendieck groups that turn 
$$ \oplusop{n\ge 0} K_0( \Bbbk [S_n] )  \ \ \  \mbox{and} \ \ \ 
   \oplusop{n\ge 0} K_0( DH_n ) $$
into graded birings, see Geissinger~\cite{Gei} for the symmetric group, 
Zelevinsky~\cite{Zel} for semisimple generalizations, and 
Bergeron and Li~\cite{BL}, Khovanov and Lauda~\cite{KL1} 
for nonsemisimple ones.  
We write \emph{birings} rather than \emph{bialgebras}, since 
these $K_0$ groups are abelian groups rather than vector spaces over some field. 
Isomorphisms in the above proposition are compatible with multiplication and 
comultiplication, and induce a biring isomorphism 
 $$  \oplusop{n\ge 0} K_0( DH_n ) \cong \oplusop{n\ge 0} K_0( \Bbbk [S_n] ) .$$ 
The second biring can be canonically identified~\cite{Gei, Zel} with an integral 
form $\Sym$ of the ring of symmetric functions in infinitely many variables. 

\vspace{0.15in} 

What can we say about $K_0$ of $\End(+^n - ^m)$ when $n,m>0$?  
Recall (end of Section~\ref{subsection-moves}) 
that $\End(+^n - ^m)$ contains 2-sided ideal $J_{n,m}$
spanned by basis elements of thickness less than $n+m$, which fits into 
short split exact sequence 
\begin{equation} \label{seq-short} 
0 \lra J_{n,m} \lra \End(+^n-^m) \lra DH_{n,m} \lra 0,  
\end{equation} 
where 
$$ DH_{n,m} \ := \ DH_n\otimes DH_m \otimes \Pi$$ 
is the tensor product of two degenerate AHA with the polynomial algebra in 
infinitely many generators. 
The earlier argument via Quillen's theorem shows that  the inclusion  
$$ \Bbbk[S_n]\otimes \Bbbk[S_m] \lra DH_n\otimes DH_m \otimes \Bbbk[c_0, 
c_1, \dots , c_r ] $$ 
induces an  isomorphisms of $K_0$-groups 
$$ K_0( \Bbbk[S_n]\otimes \Bbbk[S_m]) \cong 
 K_0( DH_n\otimes DH_m \otimes \Bbbk[c_0, c_1, \dots, c_r]).$$ 
Passing to the limit $r\to \infty$, the inclusion 
$$ \Bbbk[S_n]\otimes \Bbbk[S_m] \lra DH_{n,m}$$
induces an isomorphism of Grothendieck groups  
\begin{equation} \label{eq-K0-iso} 
  K_0( \Bbbk[S_n]\otimes \Bbbk[S_m]) \cong  K_0( DH_{n,m}).   
\end{equation} 

The split short exact sequence (\ref{seq-short}) gives rise to a split short exact sequence 
$$0 \lra K_0(J_{n,m}) \lra K_0(\End(+^n-^m)) \lra K_0(DH_{n,m}) \lra 0 $$ 
(for the definition of the $K$-group of a 2-sided ideal see~\cite{Ros, Wei})
and canonical decomposition 
$$ K_0(\End(+^n-^m)) \cong K_0(DH_{n,m}) \oplus  K_0(J_{n,m}).$$ 
Conjecture 1.1 is equivalent to the vanishing of $K_0(J_{n,m})$ for all $n,m$.

%%%%%%%%%%%%%%%%%%%%
%%%%%%%%%%%%%%%%%%%%
%%
%%   REFERENCES 
%%
%%%%%%%%%%%%%%%%%%%%
%%%%%%%%%%%%%%%%%%%%

\vspace{0.1in}
 
\noindent
{ \sl \small Mikhail Khovanov, Department of Mathematics, Columbia University, New York, NY 10027}

\noindent  {\tt \small email: khovanov@math.columbia.edu}

\end{document}